\documentclass[11pt]{article}
\usepackage{epsf}
\usepackage{epsfig}
\usepackage{amsmath}
\usepackage{amsfonts}
\usepackage{mathrsfs}
\usepackage[usenames]{color}

\newtheorem{theorem}{Theorem}

\newtheorem{condition}[theorem]{Condition}

\newtheorem{example}[theorem]{Example}

\newtheorem{lemma}[theorem]{Lemma}

\newtheorem{proposition}[theorem]{Proposition}
\newtheorem{remark}[theorem]{Remark}

\newenvironment{proof}[1][Proof]{\noindent\textbf{#1.} }{\ \rule{0.5em}{0.5em}}
\input{tcilatex}
\begin{document}

\title{Splitting for Rare Event Simulation: A Large Deviation Approach to
Design and Analysis}
\author{Thomas Dean\thanks{%
Research of this author supported in part by the National Science Foundation
(NSF-DMS-0306070 and NSF-DMS-0404806).} \,and Paul Dupuis\thanks{%
Research of this author supported in part by the National Science Foundation
(NSF-DMS-0306070 and NSF-DMS-0404806) and the Army Research Office
(W911NF-05-1-0289).} \\
Lefschetz Center for Dynamical Systems\\
Division of Applied Mathematics\\
Brown University}
\date{}
\maketitle

\begin{abstract}
Particle splitting methods are considered for the estimation of rare events.
The probability of interest is that a Markov process first enters a set $B$
before another set $A$, and it is assumed that this probability satisfies a
large deviation scaling. A notion of subsolution is defined for the related
calculus of variations problem, and two main results are proved under mild
conditions. The first is that the number of particles generated by the
algorithm grows subexponentially if and only if a certain scalar multiple of
the importance function is a subsolution. The second is that, under the same
condition, the variance of the algorithm is characterized (asymptotically)
in terms of the subsolution. The design of asymptotically optimal schemes is
discussed, and numerical examples are presented.
\end{abstract}

\section{Introduction}

The numerical estimation of probabilities of rare events is a difficult
problem. There are many potential applications in operations research and
engineering, insurance, finance, chemistry, biology, and elsewhere, and many
papers (and by now even a few books) have proposed numerical schemes for
particular settings and applications. Because the quantity of interest is
very small, standard Monte Carlo simulation requires an enormous number of
samples for the variance of the resulting estimate to be comparable to the
unknown probability. It quickly becomes unusable, and more efficient
alternatives are sought.

The two most widely considered alternatives are those based on
change-of-measure techniques and those based on branching processes. The
former is usually called \textit{importance sampling}, and the latter is
often referred to as \textit{multi-level splitting}. While good results on a
variety of problem formulations have been reported for both methods, it is
also true that both methods can produce inaccurate and misleading results.
The design issue is critical, and one can argue that the proper theoretical
tools for the design of importance sampling and splitting algorithms were
simply not available for complicated models and problem formulations.
[An alternative approach based on interacting particles has also been suggested as in \cite{delgar}.
However, we are unaware of any analysis of the performance of these schemes as the probability of interest becomes small.]

Suppose that the probability of interest takes the form $p=P\left\{ Z\in
G\right\} =\mu (G)$, where $G$ is a subset of some reasonably regular space
(e.g., a Polish space $S$) and $\mu $ a probability measure. In ordinary
Monte Carlo one generates a number of independent and identically
distributed (iid) samples $\left\{ Z_{i}\right\} $ from $\mu $, and then
estimates $p$ using the sample mean of $1_{\left\{ Z_{i}\in G\right\} }$. In
the case of importance sampling, one uses an alternative sampling
distribution $\nu $, generates iid samples $\left\{ \bar{Z}_{i}\right\} $
from $\nu $, and then estimates via the sample mean of $\left[ d\mu /d\nu %
\right] (\bar{Z}_{i})1_{\left\{ \bar{Z}_{i}\in G\right\} }$. The
Radon-Nikodim derivative $\left[ d\mu /d\nu \right] (\bar{Z}_{i})$
guarantees that the estimate is unbiased. The goal is to choose $\nu $ so
that the individual samples $\left[ d\mu /d\nu \right] (\bar{Z}%
_{i})1_{\left\{ \bar{Z}_{i}\in G\right\} }$ cluster tightly around $p$,
thereby reducing the variance. However, for complicated process models or
events $G$ the selection of a good measure $\nu $ may not be simple. The
papers \cite{glakou,glawan} show how certain standard heuristic methods
based on ideas from large deviations could produce very poor results. The
difficulty is due to points in $S$ with low probability under $\nu $ for
which $d\mu /d\nu $ is very large. The aforementioned large deviation
heuristic does not properly account for the contribution of these points to
the variance of the estimate, and it is not hard to find examples where the
corresponding importance sampling estimator is much worse than even ordinary
Monte Carlo. The estimates exhibit very inaccurate and/or unstable behavior,
though the instability may not be evident from numerical data until massive
amounts have been generated.

The most discussed application of splitting type schemes is to first
entrance probabilities, and to continue the discussion we specialize to that
case. Thus $Z$ is the sample path of a stationary stochastic process $%
\left\{ X_{i}\right\} $ (which for simplicity is taken to be Markovian), and 
$G$ is the set of trajectories that first enter a set $B$ prior to entering
a set $A$. More precisely, for disjoint $B$ and $A$ and $x\notin A\cup B$, 
\begin{equation*}
p=p(x)=P\left\{ X_{j}\in B,X_{i}\notin A,i\in \{0,\ldots ,j\},j<\infty
|X_{0}=x\right\} .
\end{equation*}%
In the most simple version of splitting, the state space is partitioned
according to certain sets $B\subset C_{0}\subset C_{1}\subset \cdots \subset
C_{K}$, with $x\notin C_{K}$ and $A\cap C_{K}=\emptyset $. These sets are
often defined as level sets of a particular function $V$, which is commonly
called an \textit{importance function}. Particles are generated and killed
off according to the following rules. A single particle is started at $x$.
Generation of particles (splitting) occurs whenever an existing particle
reaches a threshold or level $C_{i}$ for the first time. At that time, a
(possibly random) number of new particles are placed at the location of
entrance into $C_{i}$. The future evolutions of these particles are
independent of each other (and all other particles), and follow the law of $%
\{X_{i}\}$. Particles are killed if they enter $A$ before $B$. Attached to
each particle is a weight. Whenever a particle splits the weight of each
descendent equals that of the parent times a discount factor. A random tree
is thereby produced, with each leaf corresponding to a particle that has
either reached $B$ or been killed. A random variable (roughly analogous to a
single sample $\left[ d\mu /d\nu \right] (\bar{Z}_{i})1_{\left\{ \bar{Z}%
_{i}\in G\right\} }$ from the importance sampling approach) is defined as
the sum of the weights for all particles that make it to $B$. The rule that
updates the weights when a particle splits is chosen so that the expected
value of this random variable is $p$. This numerical experiment is
independently repeated a number of times, and the sample mean is again used
to estimate $p$.

There are two potential sources of poor behavior in the splitting algorithm.
The first and most troubling is that the number of particles may be large.
For example, the number could be comparable $\theta ^{K}$ for some $\theta
>1 $. In settings where a large deviation characterization of $p$ is
available, the number of levels itself usually grows with the large
deviation parameter, and so the number of particles could grow
exponentially. We will refer to this as \textit{instability} of the
algorithm. For obvious computational reasons, instability is something to be
avoided. The other source of poor behavior is analogous to that of
importance sampling (and ordinary Monte Carlo), which is high relative
variance of the estimate. If the weighting rule leads to high variation of
the weights of particles that make it to $B$, or if too many simulations
produce no particles that make it to $B$ (in which case a zero is averaged
in the sample mean), then high relative variance is likely. Note, however,
that this problem has a bounded potential for mischief, since the weights
cannot be larger than one. Such a bound does not hold for the Radon-Nikodim
derivative of importance sampling.

When the probability of interest can be approximated via large deviations,
the rate of decay is described in terms of a variational problem, such as a
calculus of variations or optimal control problem. It is well known that
problems of this sort are closely related to a family of nonlinear partial
differential equations (PDE) known as Hamilton-Jacobi-Bellman (HJB)
equations. In a pair of recent papers \cite{dupsezwan,dupwan5}, it was shown
how subsolutions of the HJB equations associated with a variety of rare
event problems could be used to construct and rigorously analyze efficient
importance sampling schemes. In fact, the subsolution property turns out to
be in some sense necessary and sufficient, in that efficient schemes can be
shown to imply the existence of an associated subsolution.

The purpose of the present paper is to show that in certain circumstances a
remarkably similar result holds for splitting algorithms. More precisely, we
will show the following under relatively mild conditions.

\begin{itemize}
\item A necessary and sufficient condition for the stability of the
splitting scheme associated to a given importance function is that a certain
scalar multiple of the importance function be a subsolution of the related
HJB equation. The multiplier is the ratio of the logarithm of the expected
number of offspring for each split and the gap between the levels.

\item If the subsolution property is satisfied, then the variance of the
splitting scheme decays exponentially with a rate defined in terms of the
value of the subsolution at a certain point.

\item As in the case of importance sampling, when a subsolution has the
maximum possible value at this point (which is the value of the
corresponding solution), the scheme is in some sense asymptotically optimal.
\end{itemize}

These results are significant for several reasons. The most obvious is that
a splitting algorithm is probably not useful if it is not stable, and the
subsolution property provides a way of checking stability. A second is that
good, suboptimal schemes can be constructed and compared via the
subsolutions framework. A third reason is that for interesting classes of
problems it is possible to construct subsolutions that correspond to
asymptotically optimal algorithms (see \cite{dupsezwan,dupwan5}).
Subsolutions can be much easier to construct than solutions. In this context
it is worth noting that the type of subsolution required for splitting (a
viscosity subsolution \cite{barcap,fleson1}) is less restrictive that the
type of subsolution required for importance sampling. Further remarks on
this point will be given in Section 5.

An outline of the paper is as follows. In the next section we describe the
probabilities to be approximated, state assumptions, and formulate the
splitting algorithm. This section also presents a closely related algorithm
that will be used in the analysis. Section 3 studies the stability problem,
and Section 4 shows how to bound the variance of an estimator in terms of
the related subsolution. The results of Sections 3 and 4 can be phrased
directly in terms of the solution to the calculus of variations problem that
is related to the large deviation asymptotics. However, for the purposes of
practical construction of importance functions the characterization via
subsolutions of a PDE is more useful. These issues are discussed in Section
5, and examples and numerical examples are presented in the concluding
Section 6.

\vspace{\baselineskip}\noindent \textbf{Acknowledgment.} Our interest in the
parallels between importance sampling and multi-level splitting was
stimulated by a talk given by P.T.\ de Boer at the RESIM conference in
Bamberg, Germany \cite{deb2}.

\section{Problem Formulation}

\subsection{Problem Setting and Large Deviation Properties}

A domain $D\subset \mathbb{R}^{d}$ is given and also a sequence of discrete
time, stationary, Markov $D-$valued processes $\left\{ X^{n}\right\} $. \
Disjoint sets $A$ and $B$ are given, and we set $\tau ^{n}\doteq \min
\left\{ i:X^{n}_i\in A\cup B\right\} $. \ The probability of interest is then%
\begin{equation*}
p^{n}(x_{n})\doteq P\left\{ X_{\tau ^{n}}^{n}\in B\left\vert
X_{0}^{n}=x_{n}\right. \right\} .
\end{equation*}%
The varying initial conditions are used for greater generality, but also
because initial conditions for the prelimit processes may be restricted to
some subset of $D$. The analogous continuous time framework can also be used
with analogous assumptions and results. For a given point $x\notin A\cup B$,
we make the following large deviation-type assumption.

\begin{condition}
\label{LDlimit}For any sequence $x_{n}\rightarrow x$,%
\begin{equation*}
\lim_{n\rightarrow \infty }-\frac{1}{n}\log p^{n}(x_{n})=W(x),
\end{equation*}%
where $W(x)$ is the solution to a control problem of the form%
\begin{equation*}
\inf \int_{0}^{t}L\left( \phi (s),\dot{\phi}(s)\right) ds.
\end{equation*}%
Here $L:\mathbb{R}^{d}\times \mathbb{R}^{d}\rightarrow \lbrack 0,\infty ]$,
and the infimum is taken over all absolutely continuous functions $\phi $
with $\phi (0)=x$, $\phi (t)\in B,\phi (s)\notin A$ for all $s\in \lbrack
0,t]$ and some $t<\infty $.
\end{condition}

\begin{remark}
\emph{The assumption that $\left\{ X^{n}\right\} $ be Markovian is not
necessary for the proofs to follow. \ For example, it could be the case that 
$X_{i}^{n} $ is the first component of a Markov process $%
(X_{i}^{n},Y_{i}^{n}) $ (e.g., so-called Markov-modulated processes). \ In
such a case it is enough that the analogous large deviation limit hold
uniformly in all possible initial conditions $Y_{0}^{n}$, and indeed the
proofs given below will carry over with only notational changes. This can be
further weakened, e.g., it is enough that the estimates hold uniformly with
sufficiently high probability in the conditioning data. \ However, the
construction of subsolutions will be more difficult, since the PDE discussed
in Section 5 is no longer available in explicit form. See \cite{dupwan5} for
further discussion on this point.}
\end{remark}

It is useful to say a few words on how one can verify conditions like
Condition \ref{LDlimit} from existing large deviation results. Similar but
slightly different assumptions will be made in various places in the sequel,
and in all cases analogous remarks will apply.

For discrete time processes one often finds process-level large deviation
properties phrased in terms of a continuous time interpolation $X^{n}(t)$,
with $X^{n}(i/n)\doteq X_{i}^{n}$ and $X^{n}(t)$ defined by piecewise linear
interpolation for $t$ not of the form $t=i/n$. In precise terms,
process-level large deviation asymptotics hold for $\left\{ X^{n}\right\} $
if the following upper and lower bounds hold for each $T\in (0,\infty )$ and
any sequence of initial conditions $x_{n}\in D$ with $x_{n}\rightarrow x$.
Define 
\begin{equation*}
I_{x}^{T}(\phi )=\int_{0}^{T}L\left( \phi (s),\dot{\phi}(s)\right) ds
\end{equation*}%
if $\phi $ is absolutely continuous with $\phi (0)=x$, and $I_{x}^{T}(\phi
)=\infty $ otherwise. If $F$ is any closed subset of $C([0,T]:D)$ then the
upper bound%
\begin{equation*}
\limsup_{n\rightarrow \infty }\frac{1}{n}\log P\left\{ X^{n}\in F\left\vert
X^{n}(0)=x_{n}\right. \right\} \leq -\inf_{\phi \in F}I_{x}^{T}(\phi )
\end{equation*}%
holds, and if $O$ is any open subset of $C([0,T]:D)$ then the lower bound%
\begin{equation*}
\liminf_{n\rightarrow \infty }\frac{1}{n}\log P\left\{ X^{n}\in O\left\vert
X^{n}(0)=x_{n}\right. \right\} \geq -\inf_{\phi \in O}I_{x}^{T}(\phi )
\end{equation*}%
holds. It is also usual to assume that for each fixed $x,T,$ and any $%
M<\infty $, the set%
\begin{equation*}
\left\{ \phi \in C([0,T]:D):I_{x}^{T}(\phi )\leq M\right\}
\end{equation*}%
is compact in $C([0,T]:D)$. The zero-cost trajectories (i.e., paths $\phi $
for which $I_{x}^{T}(\phi )=0$) are particularly significant in that all
other paths are in some sense exponentially unlikely.

With regard to Condition \ref{LDlimit}, two different types of additional
conditions beyond the sample path large deviation principle are required.
One is a condition that allows a reduction to large deviation properties
over a finite time interval. For example, suppose that there is $\bar{T}$
such that if $\phi $ enters neither $A$ nor $B$ before $\bar{T}$, then $%
I_{x}^{\bar{T}}(\phi )\geq W(x)+1$. In this case, the contribution to $%
p^{n}(x_{n})$ from sample paths that take longer than $\bar{T}$ is
negligible, and can be ignored. This allows an application of the finite
time large deviation principle. Now let $G$ be the set of trajectories that
enter $B$ at some time $t<\bar{T}$ without having previously entered $A$. By
the first condition, the asymptotic rates of decay of $p^{n}(x_{n})$ and $%
P\left\{ X^{n}\in G\left\vert X^{n}(0)=x_{n}\right. \right\} $ are the same.
The second type of condition is to impose enough regularity on the sets $A$
and $B$ and the rate function $I_{x}^{T}(\phi )$ that the infimum over the
interior and closure of $G$ are the same. These points are discussed at
length in the literature on large deviations \cite{frewen}.

\begin{example}
Assume the following conditions: $L(\cdot ,\cdot )$ is lower semicontinuous;
for each $x\in D,L(x,\cdot )$ is convex; $L(x,\cdot )$ is uniformly
superlinear; for each $x\in D$ there is a unique point $b(x)$ for which $%
L(x,b(x))=0$; $b$ is Lipschitz continuous, and all solutions to $\dot{\phi}%
=b(\phi )$ are attracted to $\theta \in A$, with $A$ open. Let $\mathcal{D}%
\subset D$ be a bounded domain that contains $A$ and $B$, and assume $%
\left\langle b(x),n(x)\right\rangle <0$ for $x\in \partial \mathcal{D}$,
where $n(x)$ is the outward normal to $\mathcal{D}$ at $x$. Suppose that the
cost to go from $x$ to any point in $\partial \mathcal{D}$ is at least $%
W(x)+1$. Then $\bar{T}$ as described above exists.
\end{example}

\subsection{The Splitting Algorithm}

In order to define a spitting algorithm we need to choose an importance
function $V(y)$ and a level size $\Delta >0$. We will require that $V(y)$ be
continuous and that $V(y)\leq 0$ for all $y\in B$. Later on we will relate $%
V $ to the value function $W$, and discuss why subsolutions to the PDE that
is satisfied by $W$ are closely related to natural candidates for the
importance function.

To simplify the presentation, we consider only splitting mechanisms with an
a priori bound $R<\infty $ on the maximum number of offspring. The
restriction is convenient for the analysis, and as we will see is without
loss of generality. The set of (deterministic) splitting mechanisms will be
indexed by $j\in \left\{ 1,\ldots ,J\right\} $. Given that mechanism $j$ has
been selected, $r(j)$ particles (with $\left\vert r(j)\right\vert \leq R$)
are generated and weights $w(j)\in {\mathbb{R}}_{+}^{r(j)}$ are assigned to
the particles. Note that we do not assume $\sum_{i=1}^{r(j)}w_{i}(j)=1$. The
class of all splitting mechanisms (i.e., including randomized mechanisms) is
identified with the set of all probability distributions on $\left\{
1,\ldots ,J\right\} $.

Associated with $V$ are the level sets 
\begin{equation*}
L_{z}=\{y\in D:V(y)\leq z\}.
\end{equation*}%
A key technical condition we use is the following. In the condition, $E_{x}$
denotes expected value given $X_{0}^{n}=x$. \ 

\begin{condition}
\label{cond:2}Let $z\in \lbrack 0,V(x)]$ be given and define $\sigma
^{n}\doteq \min \left\{ i:X^{n}_i\in A\cup L_{z}\right\} $. Then%
\begin{equation*}
\liminf_{n\rightarrow \infty }-\frac{1}{n}\log E_{x_{n}}\left[ 1_{\left\{
X_{\sigma ^{n}}^{n}\in L_{z}\right\} }\left( p^{n}(X_{\sigma
^{n}}^{n})\right) ^{2}\right] \geq W(x)+\inf_{y\in \partial L_{z}}W(y).
\end{equation*}
\end{condition}

Under the conditions discussed after Condition \ref{LDlimit} which allow one
to consider bounded time intervals, Condition \ref{cond:2} follows from the
Markov property and the large deviation upper bound.

We also define collections of sets $\left\{
C_{0}^{n}=B,C_{j}^{n}=L_{(j-1)\Delta /n},j=1,\ldots \right\} $. Define the
level function $l^{n}$ by $l^{n}(y)\doteq \min \{j\geq 0:y\in C_{j}^{n}\}$.
The location of the starting point corresponds to $l^{n}(x)=\left\lceil
nV(x)/\Delta \right\rceil $, and $l^{n}=0$ indicates entry into the target
set $B$. The splitting algorithm associated with a particular distribution $%
q $ will now be defined. Although the algorithm depends on $%
V,q,r,w,x_{n},\Delta ,A$ and $B$, to minimize notational clutter these
dependencies are not explicitly denoted.

\vspace{\baselineskip} \noindent Splitting Algorithm (SA) \vspace{%
\baselineskip}

\noindent \texttt{\indent Variables: \newline
\indent\indent$N_{r}^{n}$ number of particles in generation $r$\newline
\indent\indent$X_{r,k}^{n}$ position of $k^{th}$ particle in generation $r$%
\newline
\indent\indent$w_{r,k}^{n}$ weight of $k^{th}$ particle in generation $r$%
\newline
\indent Initialization Step: \newline
\indent\indent$N_{0}^{n}=1$, $X_{0,1}^{n}=x_{n}$, $w_{0,1}^{n}=1$\newline
\indent\indent for $r=1,\ldots ,l^{n}(x_{n})$ \newline
\indent\indent\indent$N_{r}^{n}=0$\newline
\indent\indent end\newline
\indent Main Algorithm:\newline
\indent\indent for $r=1,\ldots ,l^{n}(x_{n})$ \newline
\indent\indent\indent if $N_{r-1}^{n}=0$ then $N_{r}^{n}=0$ \newline
\indent\indent\indent else \newline
\indent\indent\indent\indent for $j=1,\ldots ,N_{r-1}^{n}$ \newline
\indent\indent\indent\indent\indent generate $Z_{r,j,i}^{n}$ a single sample
of a process with the \indent\indent\indent\indent\indent same law as $%
X_{i}^{n}$ and initial condition $Z_{r,j,0}^{n}=X_{r-1,j}^{n}$\newline
\newline
\indent\indent\indent\indent\indent let $\tau _{r,j}^{n}=\func{inf}%
\{i:Z_{r,j,i}^{n}\in A\cup C_{l^{n}(x_{n})-r}^{n}\}$ \newline
\newline
\indent\indent\indent\indent\indent Splitting Step begin\newline
\indent\indent\indent\indent\indent if $Z_{r,j,\tau_{r,j}^{n}}^{n}\notin A$%
\newline
\indent\indent\indent\indent\indent\indent let $M$ be an independent sample
from the law $q$\newline
\indent\indent\indent\indent\indent\indent for $k=1,\ldots ,|r(M)|$ \newline
\indent\indent\indent\indent\indent\indent\indent$N_{r}^{n}=N_{r}^{n}+1$ 
\newline
\indent\indent\indent\indent\indent\indent\indent$%
X_{r,N_{r}^{n}}^{n}=Z_{r,j,\tau_{r,j}^{n}}^{n}$ \newline
\indent\indent\indent\indent\indent\indent\indent$%
w_{r,N_{r}^{n}}^{n}=w_{k}(M)w_{r-1,j}^{n}$ \newline
\indent\indent\indent\indent\indent\indent end \newline
\indent\indent\indent\indent\indent end \newline
\indent\indent\indent\indent\indent Splitting Step end \newline
\newline
\indent\indent\indent\indent end \newline
\indent\indent\indent end \newline
\indent\indent end \newline
\indent Construction of a sample:\newline
\indent\indent once all the generations have been calculated we form \newline
\indent\indent the quantity \newline
\indent\indent\indent$s_{\text{SA}}^{n}=%
\sum_{j=1}^{N_{l^{n}(x_{n})}^{n}}w_{l^{n}(x_{n}),j}^{n}$. \newline
} \newline

\bigskip It should be noted that generation $1$ consists of all of the
particles that make it to set $C_{l^{n}(x)-1}^{n}$ and more generally
generation $j$ consists of all particles that make it to set $%
C_{l^{n}(x)-j}^{n}$.  We also define generation $0$ to be the initial
particle.

An estimate $\hat{p}_{\text{SA}}^{n}(x_{n})$ of $p^{n}(x_{n})$ is formed by
averaging a number of independent samples of $s_{\text{SA}}^{n}$. Observe
that once generation $r$ has been calculated the information about
generations $0$ to $r-1$ can be discarded, and so there is no need to keep
all the data in memory until completion of the algorithm. Also note that in
practice there is no need to split upon entering $C_{0}^{n}=B$.

\vspace{\baselineskip} 
\begin{figure}[htp]
\centering
\scalebox{0.8} {\input{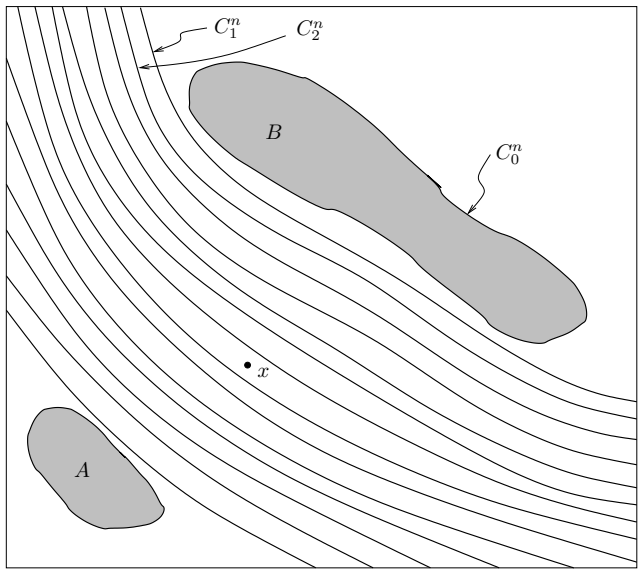}}
\caption{The Sets $A$ and $B$ and Level Sets of $V$.}
\end{figure}
\vspace{\baselineskip}

We first need to find conditions under which this splitting algorithm gives
an unbiased estimator of $p^{n}(x_n)$. To simplify this and other
calculations we introduce an auxiliary algorithm titled Splitting Algorithm
Fully Branching (SFB). The essential difference between the two is that the
process dynamics are redefined in $A$ to make it absorbing, and that
splitting continues even after a particle enters $A$. When the estimate is
constructed we only count the particles which are in $B$ in the last
generation, so that the two estimates have the same distribution. The SFB
algorithm is more convenient for purposes of analysis, because we do not
distinguish those particles which have entered $A$ from those which still
have a chance to enter $B$. Of course this algorithm would be terrible from
a practical perspective--the total number of particles is certain to grow
exponentially. However, the algorithm is used only for the purposes of
theoretical analysis, and the number of particles is not a concern. Overbars
are used to distinguish this algorithm from the previous one.

\vspace{\baselineskip} Splitting Algorithm Fully Branching (SFB) \vspace{%
\baselineskip}

\noindent \texttt{\indent Variables: \newline
\indent\indent$\bar{N}_{r}^{n}$ number of particles in generation $r$\newline
\indent\indent$\bar{X}_{r,k}^{n}$ position of $k^{th}$ particle in
generation $r$\newline
\indent\indent$\bar{w}_{r,k}^{n}$ weight of $k^{th}$ particle in generation $%
r$\newline
\indent Initialization Step: \newline
\indent\indent$\bar{N}_{0}^{n}=1$, $\bar{X}_{0,1}^{n}=x_{n}$, $\bar{w}%
_{0,1}^{n}=1$\newline
\indent\indent for $r=1,\ldots ,l^{n}(x_{n})$ \newline
\indent\indent\indent$\bar{N}_{r}^{n}=0$\newline
\indent\indent end\newline
\indent Main Algorithm:\newline
\indent\indent for $r=1,\ldots ,l^{n}(x_{n})$ \newline
\indent\indent\indent for $j=1,\ldots ,\bar{N}_{r-1}^{n}$ \newline
\indent\indent\indent\indent generate $\bar{Z}_{r,j,i}^{n}$ a single sample
of a process with the \newline
\indent\indent\indent\indent same law as $\bar{X}_{i}^{n}$ and initial
condition $\bar{Z}_{r,j,0}^{n}=\bar{X}_{r-1,j}^{n}$\newline
\newline
\indent\indent\indent\indent let $\bar{\tau}_{r,j}^{n}=\func{min}\{i:\bar{Z}%
_{r,j,i}^{n}\in A\cup C_{l^{n}(x_{n})-r}^{n}\}$ \newline
\newline
\indent\indent\indent\indent Splitting Step begin\newline
\indent\indent\indent\indent let }$\bar{M}$\texttt{\ be an independent
sample from the law }$q$\texttt{\ \newline
\indent\indent\indent\indent for $k=1,\ldots ,|r(\bar{M})|$ \newline
\indent\indent\indent\indent\indent$\bar{N}_{r}^{n}=\bar{N}_{r}^{n}+1$ 
\newline
\indent\indent\indent\indent\indent$\bar{X}_{r,\bar{N}_{r}^{n}}^{n}=\bar{Z}%
_{r,j,\bar{\tau}_{r,j}^{n}}^{n}$ \newline
\indent\indent\indent\indent\indent$\bar{w}_{r,\bar{N}_{r}^{n}}^{n}=w_{k}(%
\bar{M})\bar{w}_{r-1,j}^{n}$ \newline
\indent\indent\indent\indent end \newline
\indent\indent\indent\indent Splitting Step end \newline
\newline
\indent\indent\indent end \newline
\indent\indent end \newline
\indent Construction of a sample:\newline
\indent\indent once all the generations have been calculated we form \newline
\indent\indent the quantity \newline
\indent\indent\indent$s_{\text{SFB}}^{n}=\sum_{j=1}^{\bar{N}%
_{l^{n}(x_{n})}^{n}}1_{\left\{ \bar{X}_{l^{n}(x_{n}),j}^{n}\in B\right\} }%
\bar{w}_{l^{n}(x_{n}),j}^{n}$.\newline
}

Since the distributions of the two estimates coincide%
\begin{equation*}
E_{x_{n}}\left[ \sum_{j=1}^{N_{l^{n}(x_{n})}^{n}}w_{l^{n}(x_{n}),j}^{n}%
\right] =E_{x_{n}}\left[ \sum_{j=1}^{\bar{N}_{l^{n}(x_{n})}^{n}}1_{\left\{ 
\bar{X}_{l^{n}(x_{n}),j}^{n}\in B\right\} }\bar{w}_{l^{n}(x_{n}),j}^{n}%
\right] .
\end{equation*}%
Because of this, the SFB algorithm can be used to prove the following.

\begin{lemma}
\label{lem:unbiased}An estimator based on independent copies of $s_{\text{%
\emph{SA}}}^{n}$ is unbiased if and only if%
\begin{equation*}
E\left[ \sum_{i=1}^{r(M)}w_{i}(M)\right] =1.
\end{equation*}
\end{lemma}

\begin{proof}
It suffices to prove%
\begin{equation*}
E_{x_{n}}\left[ \sum_{j=1}^{\bar{N}_{l^{n}(x_{n})}^{n}}1_{\left\{ \bar{X}%
_{l^{n}(x_{n}),j}^{n}\in B\right\} }\bar{w}_{l^{n}(x_{n}),j}^{n}\right]
=p^{n}(x_{n}).
\end{equation*}%
We will use a particular construction of the SFB algorithm that is useful
here and elsewhere in the paper. Recall that with this algorithm every
particle splits at every generation. Hence the random number of particles
associated with each splitting can be generated prior to the generation of
any trajectories that will determine particle locations. As a consequence,
the total number of particles present at the last generation can be
calculated, as can the weight that will be assigned to each particle in this
final generation, prior to the assignment of a trajectory to the particle.
Once the weights have been assigned, the trajectories of all the particles
can be constructed in terms of random variables that are independent of
those used to construct the weights. Since the probability that any such
trajectory makes it to $B$ prior to hitting $A$ is $p^{n}(x_{n})$, 
\begin{equation*}
E_{x_{n}}\left[ \sum_{j=1}^{\bar{N}_{l^{n}(x_{n})}^{n}}1_{\left\{ \bar{X}%
_{l^{n}(x_{n}),j}^{n}\in B\right\} }\bar{w}_{l^{n}(x_{n}),j}^{n}\right]
=p^{n}(x_{n})E_{x_{n}}\left[ \sum_{j=1}^{\bar{N}_{l^{n}(x_{n})}^{n}}\bar{w}%
_{l^{n}(x_{n}),j}^{n}\right] .
\end{equation*}%
A simple proof by induction and the independence of the splitting from
particle to particle shows that%
\begin{equation*}
E_{x_{n}}\left[ \sum_{j=1}^{\bar{N}_{l^{n}(x_{n})}^{n}}\bar{w}%
_{l^{n}(x_{n}),j}^{n}\right] =\left( E\left[ \sum_{i=1}^{r(M)}w_{i}(M)\right]
\right) ^{l^{n}(x_{n})}.
\end{equation*}
\end{proof}

\bigskip

For the rest of the paper we restrict attention to splitting mechanisms that
are unbiased.

\section{Stability}

Now let an importance function $V$, level $\Delta $, and splitting mechanism 
$(q,r,w)$ be given. Define%
\begin{equation}
\mathcal{J}(x,y)\doteq \inf_{\phi ,t:\phi (0)=x,\phi
(t)=y}\int_{0}^{t}L\left( \phi (s),\dot{\phi}(s)\right) ds,
\label{TransCost}
\end{equation}%
where the infimum is over absolutely continuous functions. A function $\bar{W%
}:D\rightarrow \mathbb{R}$ will be called a \textit{subsolution} if $\bar{W}%
(x)\leq 0$ for all $x\in B$ and if $\bar{W}(x)-\bar{W}(y)\leq \mathcal{J}%
(x,y)$ for all $x,y\in D\backslash \left( A\cup B\right) $. In Section 5 we
will discuss conditions under which $\bar{W}$ can be identified as a
viscosity subsolution for an associated PDE. Recall that a splitting
algorithm is called {stable} if the total number of particles ever used
grows subexponentially as $n\rightarrow \infty $. For a given splitting
algorithm define%
\begin{equation}
\bar{W}(x)=\frac{\log Er(M)}{\Delta }V(x).  \label{ssrelation}
\end{equation}%
In this section we show that, loosely speaking, a splitting algorithm is
stable if and only if $\bar{W}$ is a subsolution.

A construction that will simplify some of the proofs is to replace a given
splitting mechanism by one for which all the weights are constant. Thus $%
(q,r,w)$ is replaced by $(q,r,\bar{w})$, where for each $j=1,\ldots ,J$ and $%
i=1,\ldots ,r(j)$,%
\begin{equation*}
\left[ \bar{w}_{i}(j)\right] ^{-1}=Er(M)=\sum_{j=1}^{J}r(j)q_{j}.
\end{equation*}%
The new splitting mechanism is also unbiased, and the distribution of the
number of particles at each stage is the same as that of $(q,r,w)$.

To establish the instability we make a very mild assumption on a large
deviation lower bound for the probability that an open ball is hit prior to
reaching $A$. This assumption can be expected to hold under conditions which
guarantee Condition \ref{LDlimit}.

\begin{proposition}
Consider an importance function $V$, level $\Delta $, and splitting
mechanism $(q,r,w)$, and define $\bar{W}$ by (\ref{ssrelation}). Suppose
there exists $y\in D\backslash \left( A\cup B\right) $ such that $\bar{W}%
(y)>0$ and%
\begin{equation}
\bar{W}(x)-\bar{W}(y)>\mathcal{J}(x,y).  \label{Instab}
\end{equation}%
Assume that $\mathcal{J}(x,y)$ is continuous at $y$. Let $\tilde{p}%
^{n}(x_{n})$ be the probability that $X^{n}$ enters the ball of radius $%
\delta >0$ about $y$ before entering $A$, given $X_{0}^{n}=x_{n}$, and
assume 
\begin{equation*}
\liminf_{n\rightarrow \infty }\frac{1}{n}\log \tilde{p}^{n}(x_{n})\geq
-\inf_{z:\left\vert z-y\right\vert <\delta }\mathcal{J}(x,z).
\end{equation*}%
Then the corresponding splitting algorithm is not stable.
\end{proposition}

\begin{proof}
It is enough to prove the instability of the algorithm that uses $(q,r,\bar{w%
})$. Since $\mathcal{J}(x,y)>0$, $V(y)<V(x)$. From the definition of $\bar{W}
$ in (\ref{ssrelation}) and (\ref{Instab}) there exist $\delta >0$ and $%
\varepsilon >0$ such that for all $z$ with $\left\vert y-z\right\vert \leq
\delta $,%
\begin{equation*}
\left[ V(x)-V(z)\right] \frac{\log Er(M)}{\Delta }>\mathcal{J}%
(x,z)+\varepsilon .
\end{equation*}%
Let $S\doteq \left\{ z:\left\vert y-z\right\vert <\delta \right\} $. By
taking $\delta >0$ smaller if necessary we can guarantee that $S\cap
A=\emptyset $ and $V(z)>0$ for all $z\in S$.

\vspace{\baselineskip} 
\begin{figure*}[htp]
\centering
\scalebox{0.8} {\input{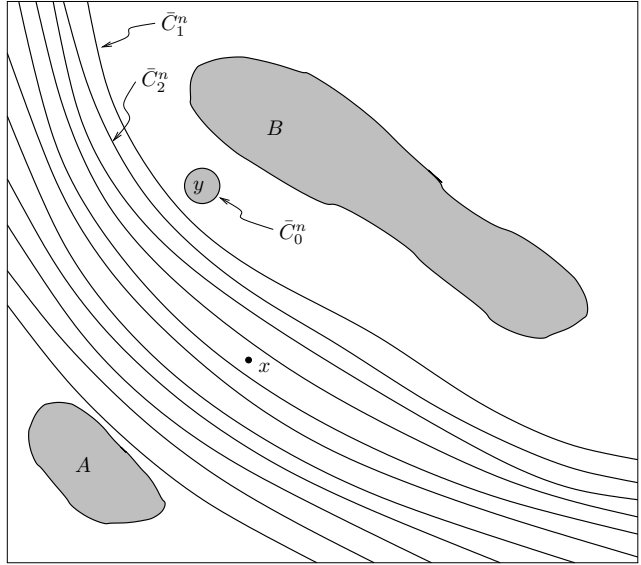}}
\caption{Level Sets of $\bar V$ in Proof of Instability.}
\end{figure*}
\vspace{\baselineskip}

Suppose one were to consider the problem of estimating $\tilde{p}^{n}(x_{n})$%
. One could use the same splitting mechanism and level sets, and even the
same random variables, except that one would stop on hitting $A$ or $S$
rather than $A$ or $B$, and the last stage would correspond to some number $%
m^{n}\leq l^{n}(x_{n})$. Of course, since $V$ is positive on $S$ it will no
longer work as an importance function, but there is a function $\bar{V}=V-a$
that will induce exactly the same level sets as $V$ and which can serve as
an importance function for this problem. See Figure 2. The two problems can
be coupled, in that exactly the same random variables can be used to
construct both the splitting mechanisms and particle trajectories up until
the particles in the $\tilde{p}^{n}(x_{n})$ problem enter $\bar{C}_{1}^{n}$.

If a particular particle has not been trapped in $A$ prior to entering $\bar{%
C}_{1}^{n}$, then that particle would also not yet be trapped in $A$ in the
corresponding scheme used to estimate $p^{n}(x_{n})$. Note also that the
number of particles that make it to $S=\bar{C}^n_0$ are at most $R$ times the number
that make it to $\bar{C}_{1}^{n}$. Let $\tilde{N}_{m^{n}}^{n}$ denote the
number of particles that make it to $S$ in the SA used to approximate $%
\tilde{p}^{n}(x_{n})$, and let $N_{l^{n}(x_{n})}^{n}$ be the number used in
the SA that approximates $p^{n}(x_{n})$. Then $N_{l^{n}(x_{n})}^{n}\geq 
\tilde{N}_{m^{n}}^{n}/R$.

Using the SFB variant in the same way that it was used in the proof of Lemma %
\ref{lem:unbiased} and that the mechanism $(q,r,\bar{w})$ is used, 
\begin{eqnarray*}
\tilde{p}^{n}(x_{n}) &=&E_{x_{n}}\left[ \sum_{j=1}^{\tilde{N}_{m^{n}}^{n}}%
\tilde{w}_{m^{n},j}^{n}\right] \\
&=&E_{x_{n}}\left[ \sum_{j=1}^{\tilde{N}_{m^{n}}^{n}}\left[ Er(M)\right]
^{-m^{n}}\right] \\
&=&\left[ Er(M)\right] ^{-m^{n}}E_{x}\left[ \tilde{N}_{m^{n}}^{n}\right] .
\end{eqnarray*}%
We now use the lower bound on $\tilde{p}^{n}(x_{n})$ and that $%
m^{n}/n\rightarrow \left[ V(x)-\sup_{z\in S}V(z)\right] /\Delta $:%
\begin{eqnarray*}
\lefteqn{\liminf_{n\rightarrow \infty }\frac{1}{n}\log E_{x_{n}}\left[ 
\tilde{N}_{m^{n}}^{n}\right] } \\
&=&\liminf_{n\rightarrow \infty }\frac{1}{n}\log E_{x_{n}}\left[ \tilde{p}%
^{n}(x_{n})\left[ Er(M)\right] ^{m^{n}}\right] \\
&\geq &-\inf_{z\in S}\mathcal{J}(x,z)+\frac{\left[ V(x)-\sup_{z\in S}V(z)%
\right] }{\Delta }\log \left[ Er(M)\right] \\
&\geq &\inf_{z\in S}\left[ \frac{\left[ V(x)-V(z)\right] }{\Delta }\log %
\left[ Er(M)\right] -\mathcal{J}(x,z)\right] \\
&\geq &\varepsilon .
\end{eqnarray*}%
It follows that 
\begin{equation*}
\liminf_{n\rightarrow \infty }\frac{1}{n}\log E_{x_{n}}\left[
N_{l^{n}(x_{n})}^{n}\right] \geq \varepsilon >0,
\end{equation*}%
which completes the proof.
\end{proof}

\bigskip

The next proposition considers stability. Here we will make a mild
assumption concerning a large deviation upper bound, which can also be
expected to hold under conditions which guarantee Condition \ref{LDlimit}.

\begin{proposition}
Consider an importance function $V$, level $\Delta $, and splitting
mechanism $(q,r,w)$, and define $\bar{W}$ by (\ref{ssrelation}). Suppose that%
\begin{equation*}
\bar{W}(x)-\bar{W}(y)\leq \mathcal{J}(x,y)
\end{equation*}%
for all $x,y\in D\backslash \left( A\cup B\right) $ and that $\bar{W}(y)\leq
0$ for all $y\in B$. Consider any $a\in [0,V(x) ]$, let $\tilde{p}%
^{n}(x_{n})$ be the probability that $X^{n}$ enters level set $L_{a}$ before
entering $A$ (given $X_{0}^{n}=x_{n}$), and assume 
\begin{equation*}
\limsup_{n\rightarrow \infty }\frac{1}{n}\log \tilde{p}^{n}(x_{n})\leq
-\inf_{z\in L_{a}}\mathcal{J}(x,z).
\end{equation*}%
Then the corresponding splitting algorithm is stable.
\end{proposition}

\begin{proof}
For each $n$ let $r^{n}$ be the value in $\left\{ 1,\ldots ,l^{n}(x)\right\} 
$ that maximizes $r\rightarrow E_{x}\left[ N_{r}^{n}\right] $. \ Since $%
r^{n}/n$ is bounded, along some subsequence (again denoted by $n$) we have $%
r^{n}/n\rightarrow v\in \lbrack 0,V(x)/\Delta ]$. Using the usual argument
by contradiction, it is enough to prove 
\begin{equation*}
\limsup_{n\rightarrow \infty }\frac{1}{n}\log E_{x_{n}}\left[ N_{r^{n}}^{n}%
\right] \leq 0
\end{equation*}%
along this subsequence. First suppose that $v=0$. Given $\delta >0$, choose $%
\bar{n}<\infty $ such that $r^{n}/n\leq \delta $ for all $n\geq \bar{n}$.
Then $N_{r^{n}}^{n,x_{n}}\leq R^{\delta n}$, and so $\limsup_{n\rightarrow
\infty }\frac{1}{n}\log E_{x_{n}}\left[ N_{r^{n}}^{n,x_{n}}\right] \leq
\delta \cdot \log R$. Since $\delta >0$ is arbitrary, this case is complete.

\vspace{\baselineskip} 
\begin{figure*}[htp]
\centering
\scalebox{0.8} {\input{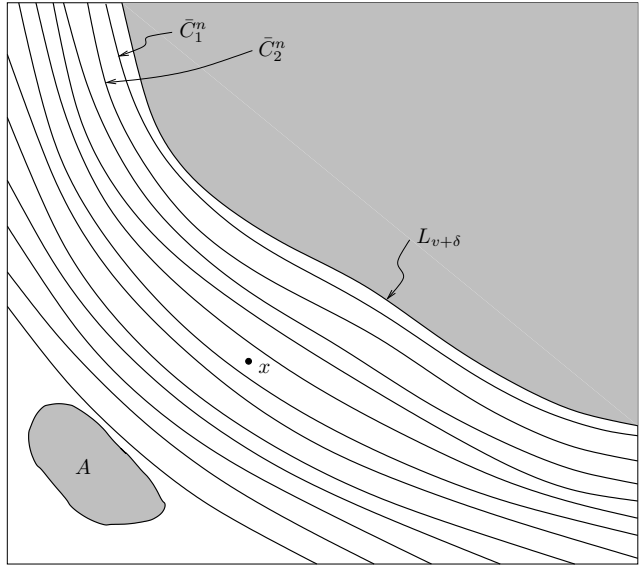}}
\caption{Level Sets of $\bar V$ in Proof of Stability.}
\end{figure*}
\vspace{\baselineskip}

Now assume $v\in (0,V(x)/\Delta ]$ and let $\delta \in (0,v)$ be given.
Suppose one were to consider the problem of estimating $\tilde{p}^{n}(x_{n})$
as defined in the statement of the proposition, with $a=V(x)-\Delta(v-\delta)$. We again use the
same splitting mechanism and level sets, except that we now stop on hitting $%
A$ or $L_{V(x)-\Delta(v-\delta)}$. An importance function with these level sets can
be found by adding a constant to $V$. We again couple the processes. \
Observe that entry into $\bar{C}_{1}^{n}$ for the $\tilde{p}^{n}(x_{n})$
problem corresponds to entry into $C_{m^{n}}^{n}$ in the ${p}^{n}(x_{n})$
problem for some $m^{n}$ such that $m^{n}/n\rightarrow [V(x)/\Delta] -(v-\delta
) $ as $n\rightarrow \infty $. Note that particles generated upon
entry into the set $C_{m^{n}}^{n}$ will correspond to generation $%
l^{n}(x)-m^{n}$ of the algorithm.  The final generation of the algorithm to
estimate  $\tilde{p}^{n}(x_{n})$ is generation $l^{n}(x)-m^{n}+1$,
corresponding to the particles generated upon reaching $\bar{C}%
_{0}^{n}=L_{V(x)-\Delta(v-\delta)}$.  Observe that since  $C_{m^{n}+1}^{n}\subset
L_{V(x)-\Delta(v-\delta)}$ every particle in the algorithm used to estimate $%
p^{n}(x_{n})$ that is not trapped in $A$ by stage  $l^{n}(x)-m^{n}+1$ will
make it to the target set $L_{V(x)-\Delta(v-\delta)}$ in the algorithm used to
estimate $\tilde{p}^{n}(x_{n})$. Hence $\tilde{N}_{l^{n}(x)-m^{n}+1}^{n}\geq
N_{l^{n}(x)-m^{n}+1}^{n}$, where $\tilde{N}_{l^{n}(x)-m^{n}+1}^{n}$
denotes the number of such particles for the SA used to estimate $\tilde{p}%
^{n}(x_{n})$.

We again use the SFB variant in the same way that it was used in the proof
of Lemma \ref{lem:unbiased} and the $(q,r,\bar{w})$ splitting mechanism to
obtain 
\begin{equation*}
\tilde{p}^{n}(x_{n})=\left[ Er(M)\right] ^{-(l^{n}(x)-m^{n}+1)}E_{x_{n}}%
\left[ \tilde{N}_{l^{n}(x)-m^{n}+1}^{n}\right] .
\end{equation*}%
It follows from $m^{n}/n\rightarrow [ V(x)/\Delta] -(v-\delta) $
that $(l^{n}(x)-m^{n}+1)/n\rightarrow \left[ v-\delta \right] $. 
Using the upper bound on $\tilde{p}^{n}(x_{n})$,%
\begin{eqnarray*}
\lefteqn{\limsup_{n\rightarrow \infty }\frac{1}{n}\log E_{x_{n}}\left[ \tilde{N}%
_{l^{n}(x)-m^{n}+1}^{n}\right] }\\ &=&\limsup_{n\rightarrow \infty }\frac{1}{n}%
\log E_{x_{n}}\left[ \tilde{p}^{n}(x_{n})\left[ Er(M)\right]
^{l^{n}(x)-m^{n}+1}\right]  \\
&\leq &-\inf_{z\in L_{_{V(x)-\Delta(v-\delta)}}}\mathcal{J}(x,z)+\left[
v-\delta \right]\log \left[ Er(M)\right]  \\
&\leq &\sup_{z\in L_{_{V(x)-\Delta(v-\delta)}}}\left[ \frac{\left[ V(x)-V(z)\right] 
}{\Delta }\log \left[ Er(M)\right] -\mathcal{J}(x,z)\right]  \\
&\leq &0.
\end{eqnarray*}%
For sufficiently large $n$ we have $r^{n}-(l^{n}(x)-m^{n}+1)\leq 2\delta n/\Delta $, and
hence $N_{r^{n}}^{n}\leq \tilde{N}_{l^{n}(x)-m^{n}+1}^{n}\cdot R^{2\delta n/\Delta }$.
It follows that 
\begin{equation*}
\limsup_{n\rightarrow \infty }\frac{1}{n}\log E_{x_{n}}\left[ N_{r^{n}}^{n}%
\right] \leq (2\delta /\Delta )\cdot \log R,
\end{equation*}%
and since $\delta >0$ is arbitrary the proof is complete.
\end{proof}

\section{Asymptotic Performance}

Since the sample $s_{\text{SA}}^{n}$ has mean $p^{n}(x_{n})$, any estimator
constructed as an average of independent copies of $s_{\text{SA}}^{n}$ is
unbiased and has variance proportional to var$_{x_{n}}\left[ s_{\text{SA}%
}^{n}\right] $. \ Once the mean is fixed, the minimization of var$_{x_{n}}%
\left[ s_{\text{SA}}^{n}\right] $ among splitting algorithms is equivalent
to the minimization of $E_{x_{n}}\left[ s_{\text{SA}}^{n}\right] ^{2}$. \ It
is of course very difficult to find the minimizer in this problem. \ When a
large deviation scaling holds, a useful alternative is to maximize the rate
of decay of the second moment, i.e., to maximize%
\begin{equation*}
\liminf_{n\rightarrow \infty }-\frac{1}{n}\log E_{x_{n}}\left[ s_{\text{SA}%
}^{n}\right] ^{2}=\liminf_{n\rightarrow \infty }-\frac{1}{n}\log E_{x_{n}}%
\left[ \sum_{j=1}^{N_{l^{n}(x_{n})}^{n}}w_{l^{n}(x_{n}),j}^{n}\right] ^{2}.
\end{equation*}%
By Jensen's inequality the best possible rate is $2W(x)$:

\begin{equation*}
\liminf_{n\rightarrow \infty }-\frac{1}{n}\log E_{x_{n}}\left[ s_{\text{SA}%
}^{n}\right] ^{2}\geq \liminf_{n\rightarrow \infty }-\frac{2}{n}\log
E_{x_{n}}\left[ s_{\text{SA}}^{n}\right] \geq 2W(x).
\end{equation*}

The main result of this section is the following.

\begin{theorem}
\label{thm:SAvar}Consider an importance function $V$, level $\Delta $, and
splitting mechanism $(q,r,w)$, and define $\bar{W}$ by (\ref{ssrelation}).
Suppose that%
\begin{equation*}
\bar{W}(x)-\bar{W}(y)\leq \mathcal{J}(x,y)
\end{equation*}%
for all $x,y\in D\backslash \left( A\cup B\right) $ and that $\bar{W}(y)\leq
0$ for all $y\in B$. \ Assume also that Conditions \ref{LDlimit} and \ref%
{cond:2} hold. Then%
\begin{equation*}
\lim_{n\rightarrow \infty }-\frac{1}{n}\log E_{x_{n}}\left[ s_{\text{\emph{SA%
}}}^{n}\right] ^{2}=W(x)-V(x)\frac{\log \left(
E\sum_{i=1}^{r(M)}w_{i}(M)^{2}\right) }{\Delta }.
\end{equation*}
\end{theorem}

\begin{proof}
It is sufficient to consider the SFB algorithm and prove that%
\begin{eqnarray*}
\lefteqn{\lim_{n\rightarrow \infty }-\frac{1}{n}\log E_{x_{n}}\left[
\sum_{j=1}^{\bar{N}_{l^{n}(x_{n})}^{n}}1_{\left\{ \bar{X}%
_{l^{n}(x_{n}),j}^{n}\in B\right\} }\bar{w}_{l^{n}(x_{n}),j}^{n}\right] ^{2}}
\\
&=&W(x)-V(x)\frac{\log \left( E\sum_{i=1}^{r(M)}w_{i}(M)^{2}\right) }{\Delta 
}.
\end{eqnarray*}%
The proof is broken into upper and lower bounds.

We first prove 
\begin{eqnarray}
\lefteqn{\limsup_{n\rightarrow \infty }-\frac{1}{n}\log E_{x_{n}}\left[
\sum_{j=1}^{\bar{N}_{l^{n}(x_{n})}^{n}}1_{\left\{ \bar{X}%
_{l^{n}(x_{n}),j}^{n}\in B\right\} }\bar{w}_{l^{n}(x_{n}),j}^{n}\right] ^{2}}
\label{upperbound} \\
&\leq &W(x)-V(x)\frac{\log \left( E\sum_{i=1}^{r(M)}w_{i}(M)^{2}\right) }{%
\Delta }.  \notag
\end{eqnarray}%
In the following display we drop cross terms to obtain the inequality, and
then use the same construction as in Lemma \ref{lem:unbiased} under which
the weights and trajectories are independent to obtain the equality. 
\begin{eqnarray*}
\lefteqn{\limsup_{n\rightarrow \infty }-\frac{1}{n}\log E_{x_{n}}\left[
\sum_{j=1}^{\bar{N}_{l^{n}(x_{n})}^{n}}1_{\left\{ \bar{X}%
_{l^{n}(x_{n}),j}^{n}\in B\right\} }\bar{w}_{l^{n}(x_{n}),j}^{n}\right] ^{2}}
\\
&\leq &\limsup_{n\rightarrow \infty }-\frac{1}{n}\log E_{x_{n}}\left[
\sum_{j=1}^{\bar{N}_{l^{n}(x_{n})}^{n}}1_{\left\{ \bar{X}%
_{l^{n}(x_{n}),j}^{n}\in B\right\} }\left( \bar{w}_{l^{n}(x_{n}),j}^{n}%
\right) ^{2}\right] \\
&=&\limsup_{n\rightarrow \infty }-\frac{1}{n}\log \left(
p^{n}(x_{n})E_{x_{n}}\left[ \sum_{j=1}^{\bar{N}_{l^{n}(x_{n})}^{n}}\left( 
\bar{w}_{l^{n}(x_{n}),j}^{n}\right) ^{2}\right] \right) .
\end{eqnarray*}%
Suppose we prove that for any $\kappa $ (and in particular $\kappa
=l^{n}(x_{n})$), that 
\begin{equation}
E_{x_{n}}\left[ \sum_{j=1}^{\bar{N}_{\kappa }^{n}}\left( \bar{w}_{\kappa
,j}^{n}\right) ^{2}\right] =\left( E\sum_{i=1}^{r(M)}w_{i}(M)^{2}\right)
^{\kappa }.  \label{secondmoment}
\end{equation}%
Since $l^{n}(x_{n})=\left\lceil nV(x_{n})/\Delta \right\rceil $, (\ref%
{upperbound}) will follow from Condition \ref{LDlimit}. \ The proof of (\ref%
{secondmoment}) is by induction. Let $M_{j}$ denote the independent random
variables used to define the splitting for the $j$th particle at stage $%
\kappa $. \ Then%
\begin{eqnarray*}
E_{x_{n}}\left[ \sum_{j=1}^{\bar{N}_{\kappa +1}^{n}}\left( \bar{w}_{\kappa
+1,j}^{n}\right) ^{2}\right] &=&E_{x_{n}}\left[ \sum_{j=1}^{\bar{N}_{\kappa
}^{n}}\left( \bar{w}_{\kappa ,j}^{n}\right)
^{2}\sum_{i=1}^{r(M_{j})}w_{i}(M_{j})^{2}\right] \\
&=&E_{x_{n}}\left[ \sum_{j=1}^{\bar{N}_{\kappa }^{n}}\left( \bar{w}_{\kappa
,j}^{n}\right) ^{2}\right] \left( E\sum_{i=1}^{r(M)}w_{i}(M)^{2}\right) \\
&=&\left( E\sum_{i=1}^{r(M)}w_{i}(M)^{2}\right) ^{\kappa +1}.
\end{eqnarray*}

We now turn to the proof of the lower bound%
\begin{eqnarray}
\lefteqn{\liminf_{n\rightarrow \infty }-\frac{1}{n}\log E_{x_{n}}\left[
\sum_{j=1}^{\bar{N}_{l^{n}(x_{n})}^{n}}1_{\left\{ \bar{X}%
_{l^{n}(x_{n}),j}^{n}\in B\right\} }\bar{w}_{l^{n}(x_{n}),j}^{n}\right] ^{2}}
\label{lowerbound} \\
&\geq &W(x)-V(x)\frac{\log \left( E\sum_{i=1}^{r(M)}w_{i}(M)^{2}\right) }{%
\Delta }.  \notag
\end{eqnarray}%
For each stage $\kappa $, let $M_{\kappa ,j}^{n}$ denote the independent
random variables used in the splitting of particle $j\in \left\{ 1,\ldots ,%
\bar{N}_{\kappa }^{n}\right\} $. \ Also, let $I_{\kappa ,j}^{n}$ denote the
disjoint decomposition of the particles in $\left\{ 1,\ldots ,\bar{N}%
_{\kappa +1}^{n}\right\} $ according to their parent particle. \ Observe
that if $k,l\in I_{\kappa ,j}^{n},k\neq l$, then for all particles descended
from $k$ and $l$, $\kappa $ is the time of their last common ancestor. \
Given $k\in I_{\kappa ,j}^{n}$, let $\bar{I}_{\kappa +1,l^{n}(x_{n}),k}^{n}$
denote the descendants of this particle at stage $l^{n}(x_{n})$. \ With this
notation we can write%
\begin{eqnarray*}
\lefteqn{E_{x_{n}}\left[ \sum_{j=1}^{\bar{N}_{l^{n}(x_{n})}^{n}}1_{\left\{ 
\bar{X}_{l^{n}(x_{n}),j}^{n}\in B\right\} }\bar{w}_{l^{n}(x_{n}),j}^{n}%
\right] ^{2}} \\
&=&\sum_{\kappa =0}^{l^{n}(x_{n})-1}E_{x_{n}}\left[ \sum_{j=1}^{\bar{N}%
_{\kappa }^{n}}\sum_{k,l\in I_{\kappa ,j}^{n},k\neq l}\sum_{m_{k}\in \bar{I}%
_{\kappa +1,l^{n}(x_{n}),k}^{n}}1_{\left\{ \bar{X}_{l^{n}(x_{n}),m_{k}}^{n}%
\in B\right\} }\bar{w}_{l^{n}(x_{n}),m_{k}}^{n}\right.  \\
&&\qquad \qquad \qquad \qquad \cdot \left. \sum_{m_{l}\in \bar{I}_{\kappa
+1,l^{n}(x_{n}),l}^{n}}1_{\left\{ \bar{X}_{l^{n}(x_{n}),m_{l}}^{n}\in
B\right\} }\bar{w}_{l^{n}(x_{n}),m_{l}}^{n}\right]  \\
&&+E_{x_{n}}\left[ \sum_{j=1}^{\bar{N}_{l^{n}(x_{n})}^{n}}1_{\left\{ \bar{X}%
_{l^{n}(x_{n}),j}^{n}\in B\right\} }\left( \bar{w}_{l^{n}(x_{n}),j}^{n}%
\right) ^{2}\right] .
\end{eqnarray*}

Let $\bar{w}_{\kappa ,j}^{n}$ denote the products of the weights accumulated
by particle $j\in \left\{ 1,\ldots ,\bar{N}_{\kappa }^{n}\right\} $ up to
stage $\kappa $, and let $\bar{w}_{\kappa +1,l^{n}(x_{n}),m}^{n}$ denote the
product of the weights accumulated by particle $m\in \left\{ 1,\ldots ,\bar{N%
}_{l^{n}(x_{n})}^{n}\right\} $ between stages $\kappa +1$ and the final
stage. Finally, let $\mathcal{F}_{\kappa }^{n}$ denote the sigma algebra
generated by $M_{s,j}^{n},s\in \left\{ 1,\ldots ,\kappa \right\} ,j\in
\left\{ 1,\ldots ,\bar{N}_{s}^{n}\right\} $ and the random variables used to
construct $\bar{X}_{s,j}^{n}$ for these same indices. Note that the future
weights are independent of $\mathcal{F}_{\kappa }^{n}$, and that the
distribution of $\bar{X}_{l^{n}(x_{n}),m}^{n}$ depends on $\mathcal{F}%
_{\kappa }^{n}$ only through $\bar{X}_{\kappa ,j}^{n}$ if $k\in I_{\kappa
,j}^{n}$ and $m\in \bar{I}_{\kappa +1,l^{n}(x_{n}),k}^{n}$. We introduce the
notation%
\begin{eqnarray*}
Y_{\kappa ,j}^{n} &\doteq &1_{\left\{ \bar{X}_{\kappa ,j}^{n}\notin
A\right\} }\left( \bar{w}_{\kappa ,j}^{n}\right) ^{2}, \\
Z_{\kappa ,k}^{n} &\doteq &\sum_{m\in \bar{I}_{\kappa
+1,l^{n}(x_{n}),k}^{n}}1_{\left\{ \bar{X}_{l^{n}(x_{n}),m}^{n}\in B\right\} }%
\bar{w}_{\kappa +1,l^{n}(x_{n}),m}^{n}
\end{eqnarray*}%
By conditioning on $\mathcal{F}_{\kappa }^{n}$ we get 
\begin{eqnarray}
\lefteqn{E_{x_{n}}\left[ \sum_{j=1}^{\bar{N}_{\kappa }^{n}}\sum_{k,l\in
I_{\kappa ,j}^{n},k\neq l}\sum_{m_{k}\in \bar{I}_{\kappa
+1,l^{n}(x_{n}),k}^{n}}1_{\left\{ \bar{X}_{l^{n}(x_{n}),m_{k}}^{n}\in
B\right\} }\bar{w}_{l^{n}(x_{n}),m_{k}}^{n}\right. }  \label{condition} \\
&&\qquad \qquad \cdot \left. \sum_{m_{l}\in \bar{I}_{\kappa
+1,l^{n}(x_{n}),l}^{n}}1_{\left\{ \bar{X}_{l^{n}(x_{n}),m_{l}}^{n}\in
B\right\} }\bar{w}_{l^{n}(x_{n}),m_{l}}^{n}\right]   \notag \\
&=&E_{x_{n}}\left[ \sum_{j=1}^{\bar{N}_{\kappa }^{n}}Y_{\kappa
,j}^{n}\sum_{k,l\in I_{\kappa ,j}^{n},k\neq l}w_{k}(M_{\kappa
,j}^{n})Z_{\kappa ,k}^{n}w_{l}(M_{\kappa ,j}^{n})Z_{\kappa ,l}^{n}\right]  
\notag \\
&=&E_{x_{n}}\left[ \sum_{j=1}^{\bar{N}_{\kappa }^{n}}Y_{\kappa
,j}^{n}\sum_{k,l\in I_{\kappa ,j}^{n},k\neq l}w_{k}(M_{\kappa
,j}^{n})w_{l}(M_{\kappa ,j}^{n})E_{\bar{X}_{\kappa ,j}^{n}}\left[ Z_{\kappa
,k}^{n}\right] E_{\bar{X}_{\kappa ,j}^{n}}\left[ Z_{\kappa ,l}^{n}\right] %
\right] .  \notag
\end{eqnarray}%
Using again the independence of the weights and trajectories as used in the
proof of Lemma \ref{lem:unbiased}, we have%
\begin{equation*}
E_{\bar{X}_{\kappa ,j}^{n}}\left[ Z_{\kappa ,k}^{n}\right] =p^{n}\left( \bar{%
X}_{\kappa ,j}^{n}\right) .
\end{equation*}%
Since 
\begin{equation*}
\mathcal{W}\doteq E\sum_{k\neq l}w_{k}(M)w_{l}(M)=E\left[ \sum_{k}w_{k}(M)%
\right] ^{2}-E\left[ \sum_{k}w_{k}(M)^{2}\right] ,
\end{equation*}%
the final expression in (\ref{condition}) equals 
\begin{equation*}
\mathcal{W}E_{x_{n}}\left[ \sum_{j=1}^{\bar{N}_{\kappa }^{n}}Y_{\kappa
,j}^{n}p^{n}\left( \bar{X}_{\kappa ,j}^{n}\right) ^{2}\right] .
\end{equation*}%
We conclude that 
\begin{eqnarray*}
\lefteqn{E_{x_{n}}\left[ \sum_{j=1}^{\bar{N}_{l^{n}(x_{n})}^{n}}1_{\left\{ 
\bar{X}_{l^{n}(x_{n}),j}^{n}\in B\right\} }\bar{w}_{l^{n}(x_{n}),j}^{n}%
\right] ^{2}} \\
&=&\mathcal{W}\sum_{\kappa =0}^{l^{n}(x_{n})-1}E_{x_{n}}\left[ \sum_{j=1}^{%
\bar{N}_{\kappa }^{n}}Y_{\kappa ,j}^{n}p^{n}\left( \bar{X}_{\kappa
,j}^{n}\right) ^{2}\right]  \\
&&+E_{x_{n}}\left[ \sum_{j=1}^{\bar{N}_{l^{n}(x_{n})}^{n}}1_{\left\{ \bar{X}%
_{l^{n}(x_{n}),j}^{n}\in B\right\} }\left( \bar{w}_{l^{n}(x_{n}),j}^{n}%
\right) ^{2}\right] .
\end{eqnarray*}

For a final time we use that the weights and trajectories are independent,
and also (\ref{secondmoment}), to argue that for any bounded and measurable
function $F$ and any stage $\kappa $,%
\begin{eqnarray*}
E_{x_{n}}\left[ \sum_{j=1}^{\bar{N}_{\kappa }^{n}}F\left( \bar{X}_{\kappa
,j}^{n}\right) \left( \bar{w}_{\kappa ,j}^{n}\right) ^{2}\right] 
&=&E_{x_{n}}\left[ F\left( \bar{X}_{\kappa ,1}^{n}\right) \right] E_{x_{n}}%
\left[ \sum_{j=1}^{\bar{N}_{\kappa }^{n}}\left( \bar{w}_{\kappa
,j}^{n}\right) ^{2}\right]  \\
&=&\left( E\left[ \sum_{i=1}^{r(M)}w_{i}(M)^{2}\right] \right) ^{\kappa
}E_{x_{n}}\left[ F\left( \bar{X}_{\kappa ,1}^{n}\right) \right] .
\end{eqnarray*}%
Thus%
\begin{eqnarray*}
\lefteqn{E_{x_{n}}\left[ \sum_{j=1}^{\bar{N}_{l^{n}(x_{n})}^{n}}1_{\left\{ 
\bar{X}_{l^{n}(x_{n}),j}^{n}\in B\right\} }\bar{w}_{l^{n}(x_{n}),j}^{n}%
\right] ^{2}} \\
&=&\mathcal{W}\sum_{\kappa =0}^{l^{n}(x_{n})-1}\left(
E\sum_{i=1}^{r(M)}w_{i}(M)^{2}\right) ^{\kappa }E_{x_{n}}\left[ 1_{\left\{ 
\bar{X}_{\kappa ,1}^{n}\notin A\right\} }p^{n}\left( \bar{X}_{\kappa
,1}^{n}\right) ^{2}\right]  \\
&&+\left( E\sum_{i=1}^{r(M)}w_{i}(M)^{2}\right) ^{l^{n}(x_{n})}E_{x_{n}}
\left[ 1_{\left\{ \bar{X}_{l^{n}(x_{n}),1}^{n}\notin A\right\} }\right] .
\end{eqnarray*}

Since $l^{n}(x_{n})$ is proportional to $n$, to prove (\ref{lowerbound}) it
is enough to show that if $\kappa _{n}$ is any sequence such that $\kappa
_{n}/n\rightarrow v\in \lbrack 0,V(x)/\Delta ]$, then%
\begin{eqnarray*}
\lefteqn{\hspace{-1in}\liminf_{n\rightarrow \infty }-\frac{1}{n}\log
E_{x_{n}}\left[ \left( E\sum_{i=1}^{r(M)}w_{i}(M)^{2}\right) ^{\kappa
_{n}}E_{x_{n}}\left[ 1_{\left\{ \bar{X}_{\kappa _{n},1}^{n}\notin A\right\}
}p^{n}\left( \bar{X}_{\kappa _{n},1}^{n}\right) ^{2}\right] \right] } \\
&\geq &W(x)-V(x)\frac{\log \left( E\sum_{i=1}^{r(M)}w_{i}(M)^{2}\right) }{%
\Delta }.
\end{eqnarray*}%
Observe that $\left\{ \bar{X}_{\kappa _{n},1}^{n}\notin A\right\} $ implies $%
\bar{X}_{\kappa _{n},1}^{n}\in C_{\left\lceil nV(x)/\Delta \right\rceil
-\kappa _{n}}^{n}$. By Condition \ref{cond:2},%
\begin{equation*}
\liminf_{n\rightarrow \infty }-\frac{1}{n}\log E_{x_{n}}\left[ 1_{\left\{ 
\bar{X}_{\kappa _{n},1}^{n}\notin A\right\} }p^{n}\left( \bar{X}_{\kappa
_{n},1}^{n}\right) ^{2}\right] \geq W(x)+\inf_{y\in \partial L_{V(x)-v\Delta
}}W(y).
\end{equation*}%
By the subsolution property, $W(y)\geq V(y)\log Er(M)/\Delta $. \ By H\"{o}%
lder's inequality 
\begin{equation*}
E\sum_{i=1}^{r(M)}w_{i}(M)^{2}\cdot Er(M)\geq E\sum_{i=1}^{r(M)}w_{i}(M)=1,
\end{equation*}%
and therefore 
\begin{equation*}
-\log \left( E\sum_{i=1}^{r(M)}w_{i}(M)^{2}\right) \leq \log \left(
Er(M)\right) .
\end{equation*}%
Thus%
\begin{eqnarray*}
\lefteqn{\liminf_{n\rightarrow \infty }-\frac{1}{n}\log E_{x_{n}}\left[
\left( E\sum_{i=1}^{r(M)}w_{i}(M)^{2}\right) ^{\kappa _{n}}E_{x_{n}}\left[
1_{\left\{ \bar{X}_{\kappa _{n},1}^{n}\notin A\right\} }p^{n}\left( \bar{X}%
_{\kappa _{n},1}^{n}\right) ^{2}\right] \right] } \\
&\geq &-v\log \left( E\sum_{i=1}^{r(M)}w_{i}(M)^{2}\right) +W(x)+\inf_{y\in
\partial L_{V(x)-v\Delta }}\frac{\log Er(M)}{\Delta }V(y) \\
&=&-v\log \left( E\sum_{i=1}^{r(M)}w_{i}(M)^{2}\right) +W(x)+\log
Er(M)\left( \frac{V(x)}{\Delta }-v\right)  \\
&\geq &W(x)-V(x)\frac{\log \left( E\sum_{i=1}^{r(M)}w_{i}(M)^{2}\right) }{%
\Delta },
\end{eqnarray*}%
and the proof is complete.
\end{proof}

\bigskip

\subsection{Design of a Splitting Algorithm}

Suppose that $V(x)$ and $\Delta $ are given and that we choose a splitting
mechanism $(q,r,w)$ which is unbiased and stable. By Theorem \ref{thm:SAvar}
the asymptotic rate of decay of the second moment is given by

\begin{equation*}
W(x)-V(x)\frac{\log \left( E\sum_{i=1}^{r(M)}w_{i}(M)^{2}\right) }{\Delta }.
\end{equation*}%
Recall from the proof of the last theorem that 
\begin{equation*}
-\log \left( E\sum_{i=1}^{r(M)}w_{i}(M)^{2}\right) \leq \log \left(
Er(M)\right) .
\end{equation*}%
\ Further equality holds if and only if $w_{i}(m)=\frac{1}{Er(M)}$ for all $%
i\in \left\{ 1,\ldots ,r(m)\right\} $ and all $m\in \left\{ 1,\ldots
,J\right\} $. \ Given the value $u=$ $Er(M)$, an alternative splitting
mechanism which is arguably the simplest which preserves the value and
achieves the equality in Holder's inequality is that defined by $J=2$ and 
\begin{equation}
q_{1}=\lceil u\rceil -u,q_{2}=1-q_{1},\;r(1)=\left\lceil u\right\rceil
,r(2)=\left\lfloor u\right\rfloor ,\;w_{i}(j)=1/u\text{ all }i,j.
\label{splitmech}
\end{equation}%
Given a subsolution $\bar{W}$, the design problem and the performance of the
resulting algorithm can be summarized as follows.

\begin{itemize}
\item Choose a level $\Delta $ and mean number of particles $u$, and define
an importance function $V$ by $\log u\cdot V(x)/\Delta =\bar{W}(x)$. \
Define the splitting mechanism by (\ref{splitmech}). \ The resulting
splitting algorithm will be stable.

\item If $s_{\text{SA}}^{n}$ is a single sample constructed according to
this algorithm, then we have the asymptotic performance%
\begin{equation*}
\lim_{n\rightarrow \infty }-\frac{1}{n}\log E_{x_{n}}[(s_{\text{SA}%
}^{n})^{2}]=W(x)+\bar{W}(x).
\end{equation*}

\item The largest possible subsolution satisfies $\bar{W}(x)=W(x)$, in which
case we achieve asymptotically optimal performance.
\end{itemize}

\begin{remark}
\label{BackingOff}\emph{Although the subsolution property guarantees
stability, it could allow for polynomial growth of the number of particles.
 If in practice one observes that a large number of particles make it to }$%
B $\emph{\ in the course of simulating a single sample }$s_{\text{SA}}^{n}$%
\emph{, then one can consider reducing the value of }$\Delta $\emph{\
slightly, while keeping the mechanism and }$V$\emph{\ fixed.  In PDE parlance,
this corresponds to the use of what is called a} strict\emph{ subsolution.  Because the value of }$\bar W(x)$\emph{ is lowered slightly, there will be a slight increase in the second moment of the estimator.  However, the strict inequality provides stronger control, and indeed the expected number of particles and moments of the number of particles will be bounded uniformly in }$n$.  
\end{remark}

\section{The Associated Hamilton-Jacobi-Bellman Equation}

The probability $p^{n}(x)$ is intimately and naturally related, via the
exponential rate $W(x)$, with a certain nonlinear PDE. This relation is well
known, and follows from the fact that $W$ is characterized in terms of an
optimal control or calculus of variations problem. We begin this section by
defining the PDE and the notion of a \textit{subsolution }in the PDE context.

Our interest in this characterization is because it is more convenient for
the explicit construction of subsolutions than the one based on the calculus
of variations problem. See, for example, the subsolutions constructed for
various large deviation problems in \cite{dupwan5} and \cite{dupsezwan}. (It
should be noted that the constructions in these papers ultimately produce 
\emph{classical} subsolutions. In contrast, the splitting algorithms require
only the weaker viscosity subsolution property. However, the smoother
subsolutions constructed in \cite{dupsezwan,dupwan5} are obtained as
mollified versions of viscosity subsolutions, and it is the construction of
these unmollified functions that is relevant to the present paper.) Other
examples will be given in the next section. Since our only interest in the
PDE is as a tool for explicit constructions, we describe the
characterization formally and in the simplest possible setting, and refer
the reader to \cite{barcap,fleson1}.

For $q\in \mathbb{R}^{d}$, let 
\begin{equation*}
\mathbb{H}(x,q)=\inf_{\beta \in \mathbb{R}^{d}}\left[ \left\langle q,\beta
\right\rangle +L(x,\beta )\right] .
\end{equation*}%
Then under regularity conditions on $L$ and the sets $A$ and $B$, $W$ can be
characterized as the maximal viscosity subsolution to 
\begin{equation*}
\mathbb{H}(x,D\bar{W}(x))=0,x\notin A\cup B,\qquad\bar{W}(x)=\left\{ 
\begin{array}{cc}
0 & x\in \partial B \\ 
\infty & x\in \partial A%
\end{array}%
\right. .
\end{equation*}%
A continuous function $\bar{W}$ is a viscosity subsolution to this equation
and boundary conditions if $\bar{W}(x)\leq 0$ for $x\in \partial B$, $\bar{W}%
(x)\leq \infty $ for $x\in \partial A$ and if the following condition holds.
If $\phi :\mathbb{R}^{d}\rightarrow \mathbb{R}$ is a smooth test function
such that the mapping $x\rightarrow \left[ \bar{W}(x)-\phi (x)\right] $
attains a maximum at $x_{0}\in \mathbb{R}^{d}\backslash \left( A\cup
B\right) $, then $\mathbb{H}(x_{0},D\phi (x_{0}))\geq 0$.
If $\bar W$ is smooth at $x_0$ then this implies $\mathbb{H}(x_{0},D\bar W (x_{0}))\geq 0$.

Note that $\mathbb{H}(x,\cdot )$ is concave for each $x_{0}\in \mathbb{R}%
^{d} $, and hence the pointwise minimum of a collection of subsolutions is
again a subsolution. It is this observation which makes the explicit
construction of subsolutions feasible in a number of interesting problems
(see \cite{dupwan5}).

In Section 2 we defined $\bar{W}(x)$ to be a subsolution to the calculus of
variations problem if it satisfied the boundary inequalities and%
\begin{equation*}
\bar{W}(x)-\bar{W}(y)\leq \mathcal{J}(x,y)
\end{equation*}%
for all $x,y\in \mathbb{R}^{d}\backslash \left( A\cup B\right) $. We now
give the elementary proof that these notions coincide. Let $\bar{W}(x)-\phi
(x)$ attain a maximum at $x_{0}$. Thus for any $\beta \in \mathbb{R}^{d}$
and all $a\in (0,1)$ sufficiently small, $\bar{W}(x_{0}+a\beta )-\phi
(x_{0}+a\beta )\leq \bar{W}(x_{0})-\phi (x_{0})$, and so 
\begin{align*}
\phi (x_{0})-\phi (x_{0}+a\beta )& \leq \bar{W}(x_{0})-\bar{W}(x_{0}+a\beta )
\\
& \leq \mathcal{J}(x_{0},x_{0}+a\beta ).
\end{align*}%
Since $\mathcal{J}(x,y)$ is defined as an infimum over all trajectories that
connect $x$ to $y$, we always have%
\begin{equation*}
\mathcal{J}(x_{0},x_{0}+a\beta )\leq \int_{0}^{a}L(x_{0}+s\beta ,\beta )ds.
\end{equation*}%
Hence if, e.g., the mapping $x\rightarrow L(x,\beta )$ is continuous, then
for all $\beta $ 
\begin{equation*}
\phi (x_{0})-\phi (x_{0}+a\beta )\leq L(x_{0},\beta )a+o(a).
\end{equation*}%
Using Taylor's Theorem to expand $\phi $, sending $a\downarrow 0$ and then
infimizing over $\beta $ gives%
\begin{equation*}
0\leq \inf_{\beta \in \mathbb{R}^{d}}\left[ \left\langle D\phi (x_{0}),\beta
\right\rangle +L(x_{0},\beta )\right] \leq \mathbb{H}(x,D\phi (x_{0})).
\end{equation*}%
Thus $\bar{W}$ is a subsolution.

The calculation just given does not show that $W$ [the solution to the calculus of variations problem] is the maximal viscosity
subsolution. Also, we have not shown that subsolutions to the PDE are also subsolutions in the calculus of variations sense. 

The characterization of $W$
as the \emph{maximal} viscosity subsolution requires that we establish $\bar{%
W}(x)\leq W(x)$ whenever $\bar{W}$ is a viscosity subsolution. \ A standard
approach to this would be to show that given a viscosity subsolution $\bar{W}
$, any point $x\notin A\cup B$, and any $\varepsilon >0$, there exists a
smooth classical subsolution $\bar{W}^{\varepsilon }$ such that $\bar{W}%
^{\varepsilon }(x)\geq \bar{W}(x)-\varepsilon $. \ When this is true the
classical verification argument \cite{fleson1} can be used to show 
$\bar W^\varepsilon (x) - \bar W^\varepsilon (y)\leq \mathcal{J}(x,y)$, and since $\varepsilon >0$ is arbitrary $\bar W$ is a subsolution to the calculus of variations problem.  This implies $W(x)\geq 
\bar{W}(x)$.

Invoking smooth subsolutions brings us very close to the method of
constructing nearly optimal importance sampling schemes as described in \cite%
{dupsezwan,dupwan5}, where the design of the scheme must be based on the
smooth classical subsolution $\bar{W}^{\varepsilon }(x)$ rather than $\bar{W}%
(x)$. In all the examples of the next subsection the inequality $\bar{W}%
(x)\leq W(x)$ can be established by constructing a nearby smooth subsolution
as in \cite{dupsezwan,dupwan5}.

\section{Numerical Examples}

In this section we present some numerical results. We study four problems:
buffer overflow for a tandem Jackson network with one shared buffer,
simultaneous buffer overflow for a tandem Jackson network with separate
buffers for each queue, some buffer overflow problems for a simple Markov
modulated queue and estimation of the sample mean of a sequence of i.i.d.
random variables.

For each case,
we present an estimate based on a stated number of runs,
standard errors and (formal) confidence intervals based on an empirical estimate of the variance,
and total computational time (intended only for comparing cases).
We also present the maximum number of particles,
which is the maximum over all runs of $\max_{r\leq l^n(x_n)}N_r$,
as well as the empirical mean and standard deviation of $\max_{r\leq l^n(x_n)}N_r$.

Subsolutions, even among those with the maximal value at a given point, are
not unique, and indeed for the problems to be discussed there are sometimes
a number of reasonable choices one could make. We will not give any details
of the proof of the subsolution property, but simply note that in each case
it can be proved by a direct verification argument as discussed in Section 5.

\subsection{Tandem Jackson Network - Single Shared Buffer}

Consider a stable tandem Jackson network with service rates $\lambda <\func{%
min}\{\mu _{1},\mu _{2}\}$. Suppose that the two queues share a single
buffer and that we are interested in the probability 
\begin{equation*}
p^{n}=P_{(0,0)}\left\{ 
\mbox{total population reaches $n$ before first
return to (0,0) }\right\}
\end{equation*}%
It is well known that 
\begin{equation*}
\lim_{n\rightarrow \infty }-\frac{1}{n}\log p^{n}=\func{min}\{\rho _{1},\rho
_{2}\}
\end{equation*}%
where $\rho _{i}=\log \frac{\mu _{i}}{\lambda }$. Further, the (continuous
time) Hamiltonian that corresponds to subsolutions of the relevant calculus
of variations problem is 
\begin{equation*}
\mathbb{H}(p)=-[\lambda (e^{-p_{1}}-1)+\mu _{1}(e^{(p_{1}-p_{2})}-1)+\mu
_{2}(e^{p_{2}}-1)].
\end{equation*}%
(see \cite{dupsezwan} for the discrete time analogue). Without loss of
generality (see \cite{dupsezwan}) one can assume that $\mu _{2}\leq \mu _{1}$%
. By inspection $\mathbb{H}(p)=0$ for $p=-\log \frac{\mu _{2}}{\lambda }%
(1,1) $ (this root is suggested by the form of the escape region), and $%
\bar W(x)=\left\langle p,x\right\rangle +\log \frac{\mu _{2}}{\lambda }$ is a
subsolution which is in fact the solution and so leads to
an asymptotically optimal splitting scheme. The table below shows the
results of a splitting simulation with 20,000 runs for $\lambda =1$, $\mu
_{1}=\mu _{2}=4.5$ and for various values of $n$. \newline

{\tiny 
\begin{tabular}{|c|c|c|c|}
\hline
$n$ & 30 & 40 & 50 \\ \hline
Theoretical Value & $2.63\times 10^{-18}$ & $1.03\times 10^{-24}$ & $%
3.80\times 10^{-31}$ \\ \hline
Estimate & $2.67\times 10^{-18}$ & $1.06\times 10^{-24}$ & $3.71\times
10^{-31}$ \\ \hline
Std. Err. & $0.11\times 10^{-18}$ & $0.05\times 10^{-24}$ & $0.20\times
10^{-31}$ \\ \hline
95\% C.I. & $[2.45,2.88]\times 10^{-18}$ & $[0.97,1.16]\times 10^{-24}$ & $%
[3.32,4.10]\times 10^{-31}$ \\ \hline
Time Taken (s) & 21 & 52 & 95 \\ \hline
Average no. particles & 23 & 31 & 37 \\ \hline
S.D. no. particles & 103 & 155 & 220 \\ \hline
Max no. particles & 2877 & 4803 & 8369 \\ \hline
\end{tabular}
}\newline

\begin{center}
Table 1. $\lambda =1$, $\mu _{1}=\mu _{2}=4.5$, asymptotically optimal
scheme.
\end{center}


It was noted in Remark \ref{BackingOff} that the number of particles
generated may grow subexponentially in $n$ and this appears to be reflected
in the data. \ Following the suggestion of the remark, we also considered a
slightly suboptimal but strict subsolution in the hopes of better controlling the
number of particles with little loss in performance. The table below shows
the results of numerical simulation for the same problem with a splitting
algorithm based on the subsolution $0.93 \times \bar W(x)$. Again
each estimate is obtained using 20,000 runs. The results are in accord with
our expectations. \newline

{\tiny 
\begin{tabular}{|c|c|c|c|}
\hline
n & 30 & 40 & 50 \\ \hline
Theoretical Value & $2.63\times 10^{-18}$ & $1.03\times 10^{-24}$ & $%
3.80\times 10^{-31}$ \\ \hline
Estimate & $2.72\times 10^{-18}$ & $1.08\times 10^{-24}$ & $3.67\times
10^{-31}$ \\ \hline
Std. Err. & $0.15\times 10^{-18}$ & $0.08\times 10^{-24}$ & $0.32\times
10^{-31}$ \\ \hline
95\% C.I. & $[2.43,3.02]\times 10^{-18}$ & $[0.93,1.23]\times 10^{-24}$ & $%
[3.04,4.30]\times 10^{-31}$ \\ \hline
Time Taken (s) & 5 & 9 & 11 \\ \hline
Average no. particles & 8 & 8 & 8 \\ \hline
S.D. no. particles & 26 & 28 & 28 \\ \hline
Max no. particles & 628 & 794 & 632 \\ \hline
\end{tabular}
\newline
}

\begin{center}
Table 2. $\lambda =1$, $\mu _{1}=\mu _{2}=4.5$, asymptotically suboptimal
scheme.
\end{center}

\subsection{Tandem Jackson Network - Separate Buffers}

\label{secGlassProb}

In the paper \cite{glaheishazaj2} the authors address the problem of
asymptotic optimality for splitting algorithms. In particular they consider
an approach to choosing level sets that are claimed to be ``consistent''
with the large deviations analysis and show that this does not always lead
to asymptotically optimal algorithms. They illustrate their results by
considering, for a tandem Jackson network, the problem of simulating the
probabilities 
\begin{equation*}
p^{n}=P_{(0,0)}\left\{ \text{both queues simultaneously exceed }n\text{
before first return to }(0,0)\right\} .
\end{equation*}%
It is shown that 
\begin{equation*}
\lim_{n\rightarrow \infty }-\frac{1}{n}\log p^{n}=\rho _{1}+\rho _{2}\doteq
\gamma ,
\end{equation*}%
and the authors propose a splitting algorithm based on the importance
function $U(x)=\gamma -\gamma \func{min}\{x_{1},x_{2}\}$, which is just a
rescaling of the target set $B = \{ (x,y) : x \geq n \mbox{ or } y \geq n \}$%
. They show that although the level sets given by this function may
intuitively seem to agree with the most likely path to the rare set
identified by the large deviations analysis, the resulting splitting
algorithm in fact has very poor performance. By analyzing this importance
function using the subsolution approach it is very easy to see why this is
the case. The Hamiltonian corresponding to subsolutions is the same as in
the previous section and it is clear to see that $U(x)$ is not a
subsolution. However, the function 
\begin{equation*}
\bar W(x)=\gamma -\rho _{1}x_{1}-\rho _{2}x_{2}
\end{equation*}%
is a subsolution. Further $\bar W(0)=\gamma $, thus the corresponding importance
function will lead to an asymptotically optimal splitting algorithm.
Numerical results are presented for the cases $\lambda =1,\mu _{1}=3,\mu
_{2}=2$ and $\lambda =1,\mu _{1}=2,\mu _{2}=3$ which are the same rates
originally considered in \cite{glaheishazaj2}. Each estimate was obtained by
a simulation using 20,000 runs.\newline

{\tiny 
\begin{tabular}{|c|c|c|c|}
\hline
$n$ & 10 & 20 & 30 \\ \hline
Theoretical Value & $9.64\times 10^{-8}$ & $1.60\times 10^{-15}$ & $%
2.64\times 10^{-23}$ \\ \hline
Estimate & $9.74\times 10^{-8}$ & $1.58\times 10^{-15}$ & $2.66\times
10^{-23}$ \\ \hline
Std. Err. & $0.17\times 10^{-8}$ & $0.03\times 10^{-15}$ & $0.06\times
10^{-23}$ \\ \hline
95\% C.I. & $[9.41,10.1]\times 10^{-8}$ & $[1.52,1.65]\times 10^{-15}$ & $%
[2.54,2.79]\times 10^{-23}$ \\ \hline
Time Taken (s) & 25 & 188 & 639 \\ \hline
Average no. particles & 25 & 53 & 81 \\ \hline
S.D. no. particles & 47 & 123 & 213 \\ \hline
Max no. particles & 550 & 1690 & 3130 \\ \hline
\end{tabular}
}

\begin{center}
Table 3. $\lambda =1,\mu _{1}=3,\mu _{2}=2$, asymptotically optimal scheme.
\end{center}

{\tiny 
\begin{tabular}{|c|c|c|c|}
\hline
$n$ & 10 & 20 & 30 \\ \hline
Theoretical Value & $9.64\times 10^{-8}$ & $1.60\times 10^{-15}$ & $%
2.64\times 10^{-23}$ \\ \hline
Estimate & $9.50\times 10^{-8}$ & $1.56\times 10^{-15}$ & $2.68\times
10^{-23}$ \\ \hline
Std. Err. & $0.26\times 10^{-8}$ & $0.06\times 10^{-15}$ & $0.13\times
10^{-23}$ \\ \hline
95\% C.I. & $[8.99,10.0]\times 10^{-8}$ & $[1.44,1.68]\times 10^{-15}$ & $%
[2.43,2.94]\times 10^{-23}$ \\ \hline
Time Taken (s) & 19 & 136 & 468 \\ \hline
Average no. particles & 26 & 54 & 86 \\ \hline
S.D. no. particles & 74 & 222 & 448 \\ \hline
Max no. particles & 1055 & 4905 & 11350 \\ \hline
\end{tabular}%
}

\begin{center}
Table 4. $\lambda =1,\mu _{1}=2,\mu _{2}=3$, asymptotically optimal scheme.
\end{center}

As expected, these results show a vast improvement over those obtained in 
\cite{glaheishazaj2}. Finally the tables below show the results of numerical
simulation for the same problem with a splitting algorithm based on the
subsolution $0.95 \times \bar W(x)$.

\vspace{\baselineskip}

{\tiny 
\begin{tabular}{|c|c|c|c|}
\hline
$n$ & 10 & 20 & 30 \\ \hline
Theoretical Value & $9.64\times 10^{-8}$ & $1.60\times 10^{-15}$ & $%
2.64\times 10^{-23}$ \\ \hline
Estimate & $9.51 \times 10^{-8}$ & $1.59\times 10^{-15}$ & $2.78\times
10^{-23}$ \\ \hline
Std. Err. & $0.20 \times 10^{-8}$ & $0.05\times 10^{-15}$ & $0.14 \times
10^{-23}$ \\ \hline
95\% C.I. & $[9.12,9.91]\times 10^{-8}$ & $[1.48,1.70]\times 10^{-15}$ & $%
[2.51,3.05]\times 10^{-23}$ \\ \hline
Time Taken (s) & 12 & 53 & 98 \\ \hline
Average no. particles & 13 & 18 & 19 \\ \hline
S.D. no. particles & 24 & 37 & 41 \\ \hline
Max no. particles & 255 & 416 & 477 \\ \hline
\end{tabular}%
}

\begin{center}
Table 5. $\lambda =1,\mu _{1}=3,\mu _{2}=2$, asymptotically suboptimal
scheme.
\end{center}

{\tiny 
\begin{tabular}{|c|c|c|c|}
\hline
$n$ & 10 & 20 & 30 \\ \hline
Theoretical Value & $9.64\times 10^{-8}$ & $1.60\times 10^{-15}$ & $%
2.64\times 10^{-23}$ \\ \hline
Estimate & $10.3 \times 10^{-8}$ & $1.52 \times 10^{-15}$ & $2.32 \times
10^{-23}$ \\ \hline
Std. Err. & $0.32 \times 10^{-8}$ & $0.08 \times 10^{-15}$ & $0.20\times
10^{-23}$ \\ \hline
95\% C.I. & $[9.66,10.9]\times 10^{-8}$ & $[1.35,1.68]\times 10^{-15}$ & $%
[1.92,2.72]\times 10^{-23}$ \\ \hline
Time Taken (s) & 10 & 37 & 67 \\ \hline
Average no. particles & 15 & 18 & 19 \\ \hline
S.D. no. particles & 39 & 61 & 72 \\ \hline
Max no. particles & 636 & 1266 & 2060 \\ \hline
\end{tabular}%
}

\begin{center}
Table 6. $\lambda =1,\mu _{1}=2,\mu _{2}=3$, asymptotically suboptimal
scheme.
\end{center}

The choice of $\bar W(x)=\gamma -\rho _{1}x_{1}-\rho _{2}x_{2}$ as a subsolution may seem arbitrary, however it turns out to be a very natural
choice. Given $\alpha > 0$ consider a ``nice'' set $B$ such that for the
subsolution $\bar W_\alpha (x) =\alpha -\rho _{1}x_{1}-\rho _{2}x_{2}$, $%
B \cap \{ x : \bar W_\alpha(x) > 0\} = \emptyset$ and $B \cap \{ x : \bar W_\alpha(x)
= 0 \} \neq \emptyset$. Then 
\begin{equation*}
\lim_{n \rightarrow \infty} - \frac{1}{n} \log p^{n} = \alpha,
\end{equation*}
where 
\begin{equation*}
p^{n} = P_{(0,0)}(\mbox{queue reaches } nB \mbox{ before first return to }
(0,0)).
\end{equation*}
Intuitively this means that all points on a level set of the function $%
\bar W_{\alpha}(x)$ have the same asymptotic probability. Thus given any such
nice set $B$ we can identify its large deviations rate by finding the unique 
$\alpha^{*}$ such that $B \cap \{ x : \bar W_{\alpha^{*}}(x) > 0\} = \emptyset$
and $B \cap \{ x : \bar W_{\alpha^{*}}(x) = 0 \} \neq \emptyset$. Further $%
\bar W_{\alpha^{*}}(x)$ will lead to an asymptotically optimal scheme.

That the family of functions $\bar W_{\alpha}$ has such a property is because the
stationary probabilities for a stable tandem Jackson network have the
product form $\pi(\{i,j\}) = (1 - \rho_{1})(1 - \rho_{2}) \rho_{1}^{i}
\rho_{2}^{j}$. Indeed, by using an argument based on the recurrence theorem,
we can see that every stable tandem Jackson network has a family of affine
subsolutions with the same property. Further this will be true for any $N$%
-dimensional queueing network for which the stationary probabilities $\pi$
have asymptotic product form, by which we mean that there exist $\rho_{1},
\ldots, \rho_{N}$ such that for any nice set $B$ 
\begin{equation*}
\lim_{n \rightarrow \infty} - \frac{1}{n} \log \pi (nB) = \inf \{ x_{1}
\rho_{1} + \cdots + x_{N} \rho_{N}: (x_{1}, \ldots, x_{N})\in B \}.
\end{equation*}

\subsection{Non-Markovian Process}

Since many models are non-Markovian we present an example of splitting for a
non-Markovian process. Consider a tandem network whose arrival and service
rates are modulated by an underlying process $M_{t}$ which takes values in
the set $\{ 1,2 \}$, such that the times taken for the modulating process to
switch states are independent exponential random variables with rate $%
\gamma(1)$ if $M$ is in state $1$ and $\gamma(2)$ otherwise. Let $%
\lambda(1), \mu_{1}(1), \mu_{2}(1)$ and $\lambda(2), \mu_{1}(2), \mu_{2}(2)$
be the service rates of the queue in the first and second states
respectively. It is known (see, e.g., \cite{dupwan4}) that the Hamiltonian
can be characterized in terms of the solution to an eigenvalue/eigenvector
problem parameterized by $p$. This characterization is used for calculating
the various roots to $\mathbb{H}(p)=0$ used below.

%
%

Consider again the single shared buffer problem. Let $\lambda(1) = 1,
\mu_{1}(1) = 3.5, \mu_{2}(1) = 2.5, \gamma(1) = 0.2$ and $\lambda(2) = 1,
\mu_{1}(2) = 4.5, \mu_{2}(2) = 4.5, \gamma(2) = 0.5$. Using a verification
argument, one can show that $\bar W(x) = 1.00029(1 - x_{1} - x_{2})$ is a
subsolution with the maximal value $\bar W(0)$. Thus using $\bar W(x)$ leads to an
asymptotically optimal splitting scheme. The results of simulations run
using this importance function are shown below, where again each estimate
was derived using 20,000 runs.

\vspace{\baselineskip} {\normalsize \vspace{0pt} }{\tiny 
\begin{tabular}{|c|c|c|c|}
\hline
n & 30 & 40 & 50 \\ \hline
Theoretical Value & $6.36 \times 10^{-13}$ & $2.88 \times 10^{-17}$ & $1.30
\times 10^{-21}$ \\ \hline
Estimate & $6.66 \times 10^{-13}$ & $2.89 \times 10^{-17}$ & $1.27 \times
10^{-21} $ \\ \hline
Std. Err. & $0.23 \times 10^{-13}$ & $0.13 \times 10^{-17}$ & $0.06 \times
10^{-21}$ \\ \hline
95\% C.I. & $[6.23,7.11]\times 10^{-13}$ & $[2.66,3.11]\times 10^{-17}$ & $%
[1.16,1.38]\times 10^{-21}$ \\ \hline
Time Taken (s) & 8 & 14 & 20 \\ \hline
Average no. particles & 4 & 5 & 5 \\ \hline
S.D. no. particles & 10 & 11 & 13 \\ \hline
Max no. particles & 188 & 214 & 280 \\ \hline
\end{tabular}
\newline
}

\begin{center}
Table 7. Markov-modulated network, total population overflow.
\end{center}

It is also worth revisiting the separate buffers problem for the same
queueing network. 
For the same arrival and service rates one can again use a verification
argument to show that 
$\bar W(x) = 2.2771 - 1.2953 x_1 - 0.9818 x_2$ leads to an asymptotically optimal
splitting scheme. Results of a simulation using 20,000 runs are shown below.

\vspace{\baselineskip} {\normalsize \vspace{0pt} }{\tiny 
\begin{tabular}{|c|c|c|c|}
\hline
n & 10 & 20 & 30 \\ \hline
Theoretical Value & $8.36 \times 10^{-10}$ & $1.07 \times 10^{-19}$ & $1.39
\times 10^{-29}$ \\ \hline
Estimate & $8.24 \times 10^{-10}$ & $1.04 \times 10^{-19}$ & $1.36 \times
10^{-29} $ \\ \hline
Std. Err. & $0.19 \times 10^{-10}$ & $0.03 \times 10^{-19}$ & $0.05 \times
10^{-29}$ \\ \hline
95\% C.I. & $[7.85,8.62]\times 10^{-10}$ & $[0.98,1.10]\times 10^{-19}$ & $%
[1.26,1.45]\times 10^{-29}$ \\ \hline
Time Taken (s) & 22 & 150 & 479 \\ \hline
Average no. particles & 31 & 59 & 89 \\ \hline
S.D. no. particles & 76 & 182 & 336 \\ \hline
Max no. particles & 1076 & 3228 & 7871 \\ \hline\hline
\end{tabular}
\newline
}

\begin{center}
Table 8. Markov-modulated network, simultaneous separate buffer overflow.
\end{center}

Finally we investigate what happens in this case if we use a strict
subsolution as importance function. The table below shows the results of a
simulation using 20,000 runs based on the importance function $
0.95 \times \bar W$.

\vspace{\baselineskip} {\normalsize \vspace{0pt} }{\tiny 
\begin{tabular}{|c|c|c|c|}
\hline
n & 10 & 20 & 30 \\ \hline
Theoretical Value & $8.36 \times 10^{-10}$ & $1.07 \times 10^{-19}$ & $1.39
\times 10^{-29}$ \\ \hline
Estimate & $8.28 \times 10^{-10}$ & $1.06 \times 10^{-19}$ & $1.43 \times
10^{-29} $ \\ \hline
Std. Err. & $0.21 \times 10^{-10}$ & $0.04 \times 10^{-19}$ & $0.07 \times
10^{-29}$ \\ \hline
95\% C.I. & $[7.87,8.69]\times 10^{-10}$ & $[0.98,1.14]\times 10^{-19}$ & $%
[1.29,1.56]\times 10^{-29}$ \\ \hline
Time Taken (s) & 15 & 70 & 157 \\ \hline
Average no. particles & 22 & 31 & 36 \\ \hline
S.D. no. particles & 50 & 88 & 116 \\ \hline
Max no. particles & 707 & 1469 & 3021 \\ \hline
\end{tabular}
\newline
}

\begin{center}
Table 9. Markov-modulated network, asymptotically suboptimal scheme.
\end{center}

\subsection{Rare Events for the Sample Mean}

It is also worth noting that this approach works just as well for finite
time problems. Assume that $X_{1},X_{2},\ldots $ is a sequence of i.i.d. N$%
(0,I^{N})$ random variables where $I^{N}$ is the $N$-dimensional identity
matrix and let $S_{n}=\frac{1}{n}\sum_{i=1}^{n}X_{i}$. Suppose that we are
interested in simulating the sequence of probabilities

\begin{equation*}
p^{n}=P\left\{ S_{n}\in C\right\}
\end{equation*}%
for some set $C$ such that $\bar{C}$ does not include the origin. For $j\in
\{1,\ldots ,n\}$ let $S_{n}(j)=\frac{1}{n}\sum_{i=1}^{j}X_{i}$. Then given
sequences $x_{n}$, $j_{n}$ and $x\in \mathbb{R}^{N}$, $t\in \lbrack 0,1]$
such that $\lim_{n \rightarrow \infty} x_{n}=x$ and $\lim_{n \rightarrow
\infty} j_{n}/n=t,$ the large deviations result 
\begin{equation*}
\lim_{n \rightarrow \infty}-\frac{1}{n}\log P\left\{ S_{n}\in
C|S_{n}(j_{n})=x_{n}\right\} =W(x,t)
\end{equation*}%
holds. Further the PDE corresponding to solutions of the calculus of
variations problem is (see \cite{dupwan5}) 
\begin{equation*}
W_{t}+\inf_{\beta }\mathbb{H}(DW;\beta )=0,
\end{equation*}%
where $\mathbb{H}(s;\beta )=\left\langle s,\beta \right\rangle +L(\beta )$
and $L(\beta )=\Vert \beta \Vert ^{2}/2$. We can put this into the general
framework in the standard way, i.e., by considering the time variable as
simply another state variable. The set $B$, for example, is then $C\times
\left\{ 1\right\} $. Strictly speaking this problem does not satisfy the
conditions used previously, since the sets $A$ and $B$ no longer have
disjoint closure. Although we omit the details, it is not difficult to work
around this problem.

It is easy to see that any affine function of the form 
\begin{equation*}
\bar{W}(x,t)=-\left\langle \alpha ,x\right\rangle +\Vert \alpha \Vert
^{2}-(1-t)H(\alpha ),
\end{equation*}%
where $H(\alpha )=\Vert \alpha \Vert ^{2}/2$, is a subsolution, though it
may not have the optimal value at $(0,0)$ and may not be less than or equal
to zero on $B$. We can use the fact that the minimum of a collection of
subsolutions is also a subsolution to build a subsolution which satisfies
the boundary condition and has the maximal value at $(0,0)$. For example,
suppose that $C=\{x\in \mathbb{R}^{2}:\left\langle p_{1},x\right\rangle \geq
1\}\cup \{x\in \mathbb{R}^{2}:\left\langle p_{2},x\right\rangle \geq 1\}$
where $p_{1}=(0.6,0.8)$ and $p_{2}=(0.6,-0.8)$. Let $\bar W_{1}(x)=1-\left\langle
p_{1},x\right\rangle -\frac{1}{2}(1-t)$, $\bar W_{2}(x)=1-\left\langle
p_{2},x\right\rangle -\frac{1}{2}(1-t)$. Then $\bar{W}=\bar W_{1}\wedge \bar W_{2}$ is
a subsolution and in fact provides an asymptotically optimal splitting
scheme since $\bar{W}(0,0)=W(0,0)$. Numerical results are shown below. Each
estimate was derived using 100,000 runs. In contrast to all the previous
examples where the process evolves on a grid, the simulated process in this
case may cross more than one splitting threshold in a single discrete time
step. This appears to increase the variance somewhat (at least if the
straightforward implementation as described in Section 2 is used), and hence
we increased the number of runs to keep the relative variances comparable.

\vspace{\baselineskip} {\normalsize \vspace{0pt} }{\tiny 
\begin{tabular}{|c|c|c|c|}
\hline
n & 20 & 30 & 40 \\ \hline
Theoretical Value & $7.75\times 10^{-6}$ & $4.33\times 10^{-8}$ & $%
2.54\times 10^{-10}$ \\ \hline
Estimate & $7.65\times 10^{-6}$ & $4.22\times 10^{-8}$ & $2.60\times
10^{-10} $ \\ \hline
Std. Err. & $0.15\times 10^{-6}$ & $0.10\times 10^{-8}$ & $0.07\times
10^{-10}$ \\ \hline
95\% C.I. & $[7.37,7.94]\times 10^{-6}$ & $[4.03,4.42]\times 10^{-8}$ & $%
[2.47,2.74]\times 10^{-10}$ \\ \hline
Time Taken (s) & 5 & 10 & 18 \\ \hline
Average no. particles & 2 & 2 & 2 \\ \hline
S.D. no. particles & 2 & 3 & 3 \\ \hline
Max no. particles & 70 & 61 & 80 \\ \hline
\end{tabular}
\newline
}

\begin{center}
Table 10. Sample mean for sums of iid.
\end{center}

%

\begin{thebibliography}{99}
\bibitem{barcap} M.~Bardi and I.~Capuzzo-Dolcetta. 
\newblock {\em Optimal Control and Viscosity Solutions of
  Hamilton-Jacobi-Bellman Equations}. \newblock Birkh{\"{a}}user, 1997.

\bibitem{deb2} P.T. de~Boer. \newblock Some observations on importance
sampling and {RESTART}. \newblock In \emph{Proc. of the 6th International
Workshop on Rare Event Simulation}, Bamberg, Germany, 2006.

\bibitem{delgar} P.~Del Moral and J.~Garnier. \newblock Genealogical particle analysis of rare events. \newblock
\emph{Ann. Appl.
Prob.},  15:2496--2534, 2005.

\bibitem{dupsezwan} P.~Dupuis, A.~Sezer, and H.~Wang. \newblock Dynamic
importance sampling for queueing networks. \newblock {\em Ann. Appl. Probab.}%
, 17:1306--1346, 2007.

\bibitem{dupwan4} P.~Dupuis and H.~Wang. \newblock Dynamic Importance
Sampling for Uniformly Recurrent {M}arkov Chains. 
\newblock {\em Ann. Appl.
Prob.}, 15:1--38, 2005.

\bibitem{dupwan5} P.~Dupuis and H.~Wang. \newblock Subsolutions of an {I}%
saacs equation and efficient schemes for importance sampling. \newblock
\emph{Math. Oper. Res.}, 32:1--35, 2007.

\bibitem{fleson1} W.~H. Fleming and H.~M. Soner. 
\newblock {\em Controlled
Markov Processes and Viscosity Solutions}. \newblock Springer-Verlag, New
York, 1992.

\bibitem{frewen} M.~I. Freidlin and A.~D. Wentzell. 
\newblock {\em Random
Perturbations of Dynamical Systems}. \newblock Springer-Verlag, New York,
1984.

\bibitem{glaheishazaj2} P.~Glasserman, P.~Heidelberger, P.~Shahabuddin, and
T.~Zajic. \newblock A large deviations perspective on the efficiency of
multilevel splitting. \newblock {\em IEEE Trans. Automat. Control},
43:1666--1679, 1998.

\bibitem{glakou} P.~Glasserman and S.~Kou. \newblock Analysis of an
importance sampling estimator for tandem queues. 
\newblock {\em ACM Trans.
Modeling Comp. Simulation}, 4:22--42, 1995.

\bibitem{glawan} P.~Glasserman and Y.~Wang. \newblock Counter examples in
importance sampling for large deviations probabilities. 
\newblock {\em Ann.
Appl. Prob.}, 7:731--746, 1997.
\end{thebibliography}

\end{document}